\documentclass[AMA,STIX1COL,12pt]{WileyNJD-v2}

\articletype{Article Type}%

\received{}
\revised{}
\accepted{}

\hypersetup{hidelinks}
\usepackage{todonotes}
\usepackage{siunitx}
\usepackage{changes}
\usepackage{lineno}
\DeclareMathOperator{\tr}{tr}
\DeclareMathOperator{\argmin}{arg\;min}
\begin{document}

\title{An asynchronous variational integrator for the phase field approach to dynamic fracture}

\author[1]{Zongwu Niu}

\author[1,3]{Vahid Ziaei-Rad}

\author[1]{Zongyuan Wu}

\author[1,2]{Yongxing Shen*}

\authormark{Zongwu Niu \textsc{et al.}}

\address[1]{\orgdiv{University of Michigan -- Shanghai Jiao Tong University Joint Institute}, \orgname{Shanghai Jiao Tong University}, \orgaddress{\state{Shanghai}, \country{China}}}

\address[2]{\orgname{Shanghai Key Laboratory for Digital
Maintenance of Buildings and Infrastructure}, \orgaddress{\state{Shanghai}, \country{China}}}

\address[3]{\orgdiv{Department of Environmental Informatics}, \orgname{Helmholtz Centre for Environmental Research -- UFZ}, \orgaddress{\state{Leipzig}, \country{Germany}}}

\corres{*Yongxing Shen, University of Michigan -- Shanghai Jiao Tong University Joint Institute, Shanghai Jiao Tong University, Shanghai 200240, China. \email{yongxing.shen@sjtu.edu.cn}}


\abstract[Summary]{The phase field approach is widely used to model fracture behaviors due to the absence of the need to track the crack topology and the ability to predict crack nucleation and branching. In this work, the asynchronous variational integrators (AVI) is adapted for the phase field approach of dynamic brittle fractures. The AVI is derived from Hamilton's principle and allows each element in the mesh to have its own local time step that may be different from others'. While the displacement field is explicitly updated, the phase field is implicitly solved, with upper and lower bounds strictly and conveniently enforced. In particular, the AT1 and AT2 variants are equally easily implemented. Three benchmark problems are used to study the performances of both AT1 and AT2 models and the results show that AVI for phase field approach significantly speeds up the computational efficiency and successfully captures the complicated dynamic fracture behavior.}



\maketitle



\section{Introduction}\label{sec1}
\linenumbers
Dynamic fracture refers to  crack development processes accompanied by fast changes in applied loads and rapid crack propagation, where inertial forces play an important role during the evolution. Application examples of dynamic fracture include drop tests of electronic devices, oil exploitation,  and impact on vehicles. 

Dynamic fracture of solids has been extensively studied \cite{freund1998dynamic,ravi2004dynamic,fineberg2015recent,SUN2021Astate}. 
Over the past decades, various numerical methods\cite{rabczuk2013computational}  to simulate dynamic fracture have been proposed, which can be classified into two groups: discrete approaches and smeared-crack ones. A discrete approach explicitly describes the crack topology, such as the extended finite element method \cite{reth2005An}, cohesive zone model \cite{Nguyen2014Discontinuous}, element deletion method \cite{Liu2014ARegular}, cracking element method\cite{ZHANG2019cracking}, and phantom nodes method\cite{song2006Amethod}, just to name a few in the context of dynamic fracture. Conversely, a smeared-crack approach represents the crack by a smeared crack band, which includes the gradient damage model \cite{Li2016Gradient}, the thick level set approach \cite{Moreau2015Explicit}, and so on.

The regularized variational fracture method \cite{BOURDIN2000797}, also called the phase field method for fracture, belongs to the group of smeared-crack approaches. It originates from Griffith's energetic theory and was developed based on the variational approach to brittle fracture by Francfort and Marigo \cite{FRANCFORT1998revisit}. The formulation solves crack problems by minimizing an energy functional that consists of the elastic energy, external work, and the crack surface energy. This way, crack evolution is a natural outcome of the solution. 
The phase field method possesses the following advantages:  (1) the crack evolves naturally and there is no need of a crack tracking algorithm; (2) there is no need of additional criterion for crack branching and merging; (3) the implementation is straightforward even for complicated crack problems in 3D. These advantages facilitate its application for various fracture problems, such as shell fracture\cite{AMIRI2014thinshell}, beam fracture\cite{LAI2020phase}, and carbon dioxide  fracturing\cite{MOLLAALI2019Numerical}. For  details about the implementation of this approach, we refer the reader to the work by Shen {et al.}~\cite{shen2018implementation} 


For dynamic fracture, Borden et al.~\cite{BORDEN201277} combined the phase field method with isogeometric analysis with local adaptive refinement to simulate dynamic brittle fracture. Nguyen and Wu~\cite{Vinh2018modeling} presented a phase field regularized cohesive zone model for dynamic brittle fracture. Hao et al.~\cite{HAO2022phasefield,HAO2022aphase} developed formulations for high-speed impact problems for metals, accounting for volumetric and shear fractures. 

However, the phase field model suffers from high computational cost partly because of the small critical time step, which, in turn, results from the necessary fine spatial discretization near the crack to resolve the regularization length scale. In order to overcome this challenge, various schemes to accelerate such computation have been proposed. Tian et al.~\cite{Tian2019hybird} presented a multilevel hybrid adaptive finite element phase field method for quasi-static and dynamic brittle fracture, wherein the refinement is based on the crack tip identified with a certain scheme. Ziaei-Rad and Shen\cite{ZIAEIRAD2016Massive} developed a parallel algorithm on the graphical processing unit with time adaptivity strategy to speed up the computation.  
Li et al.~\cite{Li2019variational} proposed a variational h-adaption method with both a mesh refinement and a coarsening scheme based on an energy criterion. 
Engwer et al.~\cite{Thomas2021Lscheme} proposed a linearized staggered scheme with dynamic adjustments of the stabilization parameters throughout the iteration to reduce the computational cost. 


In this work, with the aim of accelerating dynamic phase field fracture computations, we adapt the asynchronous variational integrator (AVI) for this problem. The AVI is an instance of variational integrators. Variational integrators are a class of time integration algorithms derived from Hamilton's principle of stationary action and have the advantages of symplectic-momentum conservation and remarkable energy (or Hamiltonian) behavior for long-time integration. In essence, they can be classified into synchronous variational integrators and asynchronous variational integrators.  The former, such as central difference, requires all unknown variables to be solved with the same time step, taking into account global requirement of stability and accuracy. 

In contrast, the latter allows independent time steps for each term contributing to the action functional, effectively independent time steps for each element in the context of finite elements. This asynchrony allows the elements with smaller time steps to be more frequently updated. Moreover, the method may be made fully explicit and even in the implicit case,   only assembly of local reaction forces and stiffness matrices instead of global ones is needed. For linear elastodynamics, the AVI was first introduced by Lew et al.~\cite{lew2003asynchronous,lew2004variational}, and the stability and convergence of AVI have been proved by Fong et al.~\cite{FONG2008Stability} and by Focardi and Mariano \cite{Matteo2008Convergence}, respectively. In addition, the AVI has been extended to the contact problem\cite{Ryckman2012AVIcontact}, wave propagation\cite{Liu2020Extended} and computer graphics\cite{Thomas2008asyncloth}. 


In the case of AVI for the phase field approach to fracture, a few adjustments need to be made. First and foremost, the overall Lagrangian is free of the time derivative of the phase field; hence solving the phase field is a local steady-state problem. More precisely, the coupled multi-field system is solved by employing a  staggered scheme, in which the displacement and velocity fields are integrated with an explicit scheme while the phase field is the solution of an inequality-constrained optimization problem.  In essence, the phase field of only one element is solved at a time, for which it is very convenient to enforce the inequality constraint compared to doing so for the entire domain. This feature allows implementing the AT1 variant of the method \cite{Pham2011Gradient} with a similar cost for the more widely used AT2 variant.

There are other asynchronous methods for dynamic fracture with the phase field. For example, Ren et al.~\cite{REN2019explict} proposed an explicit phase field formulation where the mechanical field is solved with a larger time step while the phase field is updated with smaller sub-steps. Suh and Sun\cite{SUH2021asynchronous} presented a subcycling method to capture the brittle fracture in porous media, where the heat transfer between the fluid and solid constituents is solved with different time steps as integer multiples of each other.  Note that the methods in Refs.~\citenum{REN2019explict,SUH2021asynchronous} are not variational, hence may not enjoy the said advantages of variational integrators.

The paper is organized as follows. In Section \ref{sec2}, we briefly review the formulation of the phase field model for brittle fracture and introduce  Hamilton's principle in the continuum Lagrangian framework. In Section \ref{sec3}, we present the asynchronous spacetime discretization scheme and derive discrete Euler-Lagrange equations by discrete variational principle, then we present how to solve the mechanical field and phase field by a staggered scheme. In addition, we  summarize the overall implementation of AVI for the phase field approach to fracture. In Section \ref{sec4}, we showcase three benchmark examples under dynamic loading for verification  and examining the performance. Finally, we conclude this work in Section \ref{sec5}.

\section{Formulations }\label{sec2}
This section devotes to the formulation of a dynamic fracture phase field model  through Hamilton's principle for an elastic body with possible cracks represented by a phase field.


\subsection{Hamilton's principle}

Let $\Omega\subset\mathbb{R}^n$, $n=2,3$, be the domain occupied by the reference configuration of a body with possible cracks. Hamilton's principle states that the true trajectory of a body with prescribed initial and final conditions is the stationary point of the action functional with respect to arbitrary admissible variations. Here, we consider  $\Omega $ with possible internal cracks during a specified time interval $ t\in [t_0, t_f]$  with the action functional given by
\begin{equation}
    S(\boldsymbol{u},d) = \int_{t_0}^{t_f} L(\boldsymbol{u},\boldsymbol{\dot u}, d) \,\mathrm{d}t,
\end{equation}
where $ \boldsymbol{u}(\boldsymbol{X},t) $, $\boldsymbol{X}\in\Omega$, denotes the displacement field of the body, and $\boldsymbol{\dot u}=\mathrm{d}\boldsymbol{u}/\mathrm{d}t$ is the velocity field. The scalar field $d:\Omega\times[t_0, t_f]\rightarrow[0,1]$ is called the phase field, which approximates possible sharp cracks in a diffusive way. Herein, the Lagrangian function is in the form
\begin{equation}\label{eq:Lagrangian}
     L(\boldsymbol{u},\boldsymbol{\dot u}, d)  = T(\boldsymbol{\dot u})- V(\boldsymbol{u}, d) - \Gamma(d),
\end{equation}
where $V(\boldsymbol{u}, d) $ is the potential energy,  $\Gamma(d)$ is the crack surface energy, and
\begin{equation}
    T(\boldsymbol{\dot u}) = \int_{\Omega} \frac{1}{2}\rho \boldsymbol{\dot u} \cdot \boldsymbol{\dot u}\,\mathrm{d}\Omega
    \label{eq:kinetic}
\end{equation}
is the kinetic energy, where $\rho$ is the initial mass density.

\subsection{Phase field approximation}
\begin{figure}[htbp] 
\centerline{\includegraphics[width=0.6\textwidth]{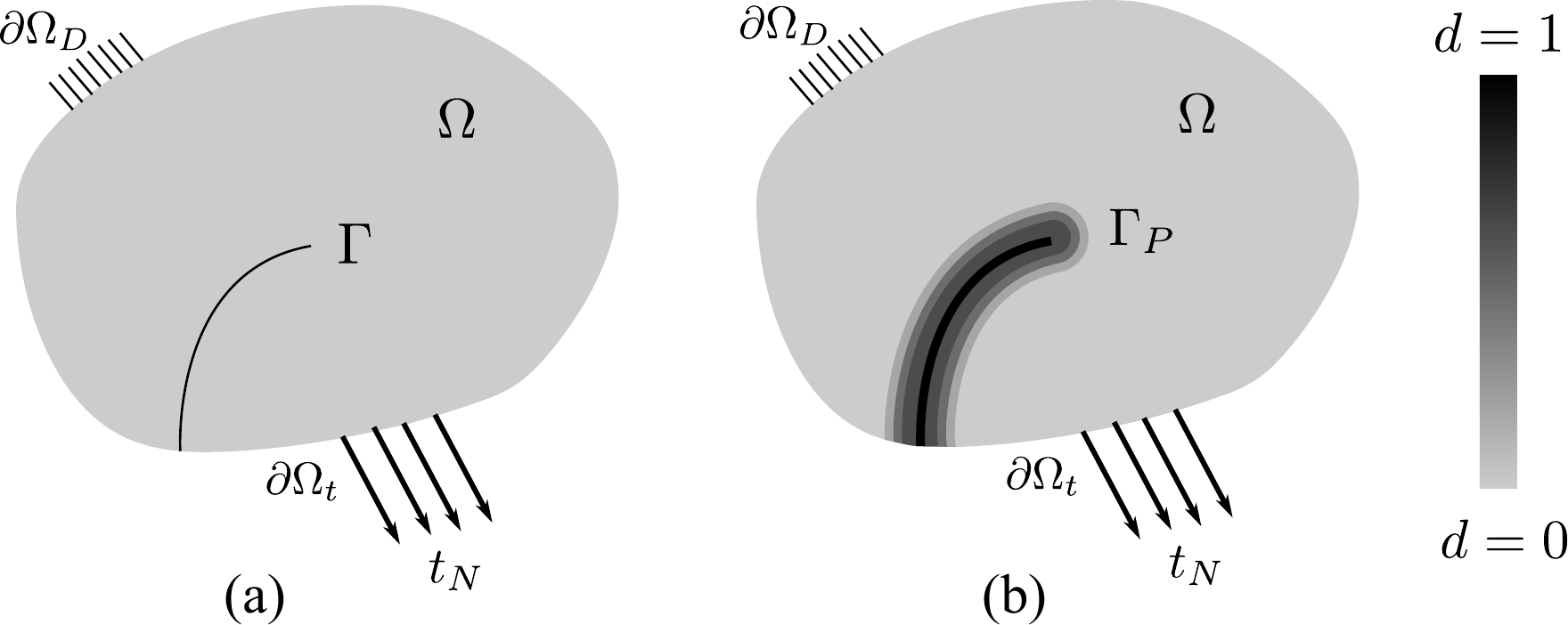}}
\caption{Body with an internal crack with: (a) a sharp crack; (b) a  crack approximated by the phase field.
\label{fig: phase field}}
\end{figure}

In this subsection, we revisit the two versions of the phase field model as a basis for subsequent development. In the phase field model of fracture, the sharp crack surface is approximated by  a scalar phase field $ d $ as shown in Figure \ref{fig: phase field}. The range of this field $d$ has to be between 0 and 1. In particular, our convention is such that the region with $d =1 $ represents the fully cracked state and that with $ d =0 $ represents the pristine state of the material. Following Ref.~\citenum{BOURDIN2000797},  the crack surface energy $\Gamma(d)$ in \eqref{eq:Lagrangian} is given by
\begin{equation}\label{eq:Gamma_d}
   \Gamma(d) = \int_{\Omega} g_c \gamma (d, \nabla d) ~\mathrm{d}\Omega, 
\end{equation}
where $ g_c > 0 $ is the critical crack energy release rate, and $\gamma(d,\nabla d)$ is the crack surface density per unit volume,
\begin{equation}\label{eq:crack density}
     \gamma(d,\nabla d) = 
     \frac{1}{4c_w}\left( \frac{w(d)}{\ell}+\ell |\nabla d|^2 \right),
\end{equation}
where $ \ell > 0 $ is the regularization length scale parameter, which controls the width of the transition region of the smoothed crack. Crack geometric function $w(d)$ and normalization constant $c_w = \int_0^1 \sqrt{w(d)}\,\mathrm{d}d $ are model dependent. Specifically, for the brittle fracture, classical examples are $w(d)=d^2$ and $c_w = 1/2$  for the AT2 model; and $w(d)=d$ and  $c_w = 2/3$  for the AT1 model \cite{Pham2011Gradient}. In addition, a notable difference between the AT2 and AT1 models is that the former gives rise to a more diffuse phase field profile while the latter generates a phase field profile with a narrow {support} near the crack.
 
The  potential energy $V(\boldsymbol{u},d)$ in \eqref{eq:Lagrangian} is expressed as
\begin{equation} 
    V(\boldsymbol{u},d) = \int_{\Omega} \psi(\boldsymbol{\varepsilon}(\boldsymbol{u}),d)\,\mathrm{d} \Omega - \int_{\partial\Omega_t } \boldsymbol{t}_{N} \cdot \boldsymbol{u} \,\mathrm{d}A - \int_{\Omega} \boldsymbol{b} \cdot \boldsymbol{u}\,\mathrm{d} \Omega,
    \label{eq:regular potential}
\end{equation}
where $\psi$ is the strain energy density,   $\boldsymbol{t}_{N}:\partial\Omega_t\rightarrow\mathbb{R}^n$ is the prescribed traction boundary condition{,} and $ \boldsymbol{b}:\Omega\rightarrow\mathbb{R}^n$ is the body force. The strain tensor is given by $\boldsymbol{\varepsilon} = ( \nabla \boldsymbol{u} + \nabla  \boldsymbol{u}^T) /2 $, where $\nabla(\cdot)$ is the gradient operator with respect to $\boldsymbol{X}$.

Here we adopt a form for $\psi$ that accounts for the unilateral constraint following Miehe et al.~\cite{MIEHE20102765} which involves spectral decomposition of $\boldsymbol{\varepsilon}$. Other choices are, for example, the  volumetric-deviatoric split  by Amor et al.~\cite{AMOR20091209}, the micromechanics-informed model by Liu et al.~\cite{LIU2021107358}{,} and the model by Wu et al.~\cite{WU2020112629} In the chosen formulation, the strain energy density takes the following form
\begin{equation}
    \psi(\boldsymbol{\varepsilon},d ) = g(d) \psi_{+} + \psi_{-}, 
\end{equation} 
where $g(d) = (1-d)^2 $ is the degradation function, and $\psi_+$ and $\psi_-$ are, respectively, the crack-driving and persistent portions of the strain energy density as
\begin{equation}
   \psi_{\pm}(\boldsymbol{\varepsilon}) = \frac{\lambda}{2} \langle \tr (\boldsymbol{\varepsilon}) \rangle_{\pm}^2 + \mu\tr\left(\boldsymbol{\varepsilon}^2_{\pm}\right),
\end{equation}
where $\lambda$ and $\mu$ are Lam\'e constants such that $\mu>0$ and $\lambda + 2\mu >0$, the Macauley  bracket is defined as  $ \langle \cdot \rangle _{\pm} = (\cdot \pm \lvert \cdot \rvert )/2  $, and
\begin{equation}
   \boldsymbol{\varepsilon}_{\pm} = \sum\limits_{a=1}^{3} \langle \varepsilon_a \rangle _{\pm} \mathbf{n}_a \otimes \mathbf{n}_a,
\end{equation}
where  $\{\varepsilon_a\}_{a=1}^3$ denote the  principal strains, $\mathbf{n}_a$ are the corresponding orthonormal principal directions, and the operator $\otimes$ represents the dyadic product.  Correspondingly, the Cauchy stress tensor is 
\begin{equation}
    \boldsymbol{\sigma}_{\pm}(\boldsymbol{\varepsilon})=\frac{\partial\psi_{\pm}}{\partial\boldsymbol{\varepsilon}}= \lambda \langle \tr (\boldsymbol{\varepsilon}) \rangle _{\pm} \boldsymbol{1} + 2\mu \sum\limits_{a=1}^{3}\langle \varepsilon_a \rangle_{\pm} \mathbf{n}_a \otimes \mathbf{n}_a,
\end{equation}
where $\boldsymbol{1} $ is the second-order identity tensor.

\subsection{Spatial discretization}
In this subsection, we obtain the semi-discrete Lagrangian by discretizing the displacement field and the phase field with a finite element mesh $\mathcal{T}_h$ for $\Omega$. Let $\eta$ be the set of nodes of $\mathcal{T}_h$.  The discretized fields take the following form
\begin{equation}\label{Eq:phi_d}
	     \boldsymbol{u}(\boldsymbol{X})=\sum_{a\in\eta} N_a(\boldsymbol{X}) \mathbf{u}_a,\quad
	     d(\boldsymbol{X})=\sum_{a\in\eta}  N_a(\boldsymbol{X}) d_a,
\end{equation} 
where $\mathbf{u}_a \in\mathbb{R}^n$ and $d_a\in\mathbb{R}$ are the displacement vector and phase field {value} at node $a \in \eta$, respectively, and $N_a$ is the finite element shape function associated with node $a$. In this work, quadrilateral finite elements and {the} associated first-order shape functions are used. 

The Lagrangian $L$  may be decomposed as
\begin{equation}\label{Eq:L_e}
		L(\boldsymbol{u},\boldsymbol{\dot u}, d)=\sum_{e\in\mathcal{T}_h}L_e\left(\mathbf{u}_e,\dot{\mathbf{u}}_e,\mathbf{d}_e\right),
\end{equation}
where $e$ is an element of the mesh $\mathcal{T}_h$, and $\mathbf{u}_e$, $\dot{\mathbf{u}}_e$, and $\mathbf{d}_e$ are the vectors containing the displacements, velocities, and phase field values of all the nodes of element $e$, respectively. The quantity $L_e$ is given by 
\begin{equation}\label{Eq:Lagrangian_loc}
		L_e\left(\mathbf{u}_e,\dot{\mathbf{u}}_e,\mathbf{d}_e\right)=T_e\left(\dot{\mathbf{u}}_e\right)-V_e\left(\mathbf{u}_e,\mathbf{d}_e\right) - \Gamma_e (\mathbf{d}_e),
\end{equation}
where $V_e$ and $ \Gamma_e$ are the elemental potential energy and elemental surface energy, respectively, and
\begin{equation}\label{Eq:Kinetic_loc}
		T_e\left(\dot{\mathbf{u}}_e\right)=\frac12\dot{\mathbf{u}}^T_e\mathbf{m}_e\dot{\mathbf{u}}_e
\end{equation}
is the elemental kinetic energy, where $\mathbf{m}_e$ is the diagonal element mass matrix. Hence, the space-discretized action is in the form
\begin{equation}\label{eq: spacedis_action}
    S\left( \{\mathbf{u}_e\},\{\dot{\mathbf{u}}_e\},\{\mathbf{d}_e\} \right) = \int_{t_0}^{t_f} \sum_{e\in\mathcal{T}_h} \Big( \frac12\dot{\mathbf{u}}^T_e\mathbf{m}_e\dot{\mathbf{u}}_e - V_e\left(\mathbf{u}_e,\mathbf{d}_e\right) - \Gamma_e (\mathbf{d}_e)\Big)\,\mathrm{d}t.
\end{equation}

For later convenience, we let $\eta(e)$ denote the set of nodes of $e$ and $\eta^{-1}(a)$ the set of elements that contain node $a$.

\section{Asynchronous variational integrator with the fracture phase field}\label{sec3}

The main feature of the AVI is to assign different time steps to different elements of $\mathcal{T}_h$. The key idea is the stationarity of \eqref{eq: spacedis_action}, a functional over space and time, which allows to divide the total Lagrangian into contributions from elemental terms which may possess independent time steps. For example, the smaller elements in the mesh may be updated a few times while the larger elements are held, according to either a preset schedule or a schedule determined on the fly.

In the context of fracture phase field in this work, the phase field is implicitly solved at the element level instead of at the global level, permitting more efficient solvers of inequality constraints.

In this section, we detail a phase field formulation and implementation for dynamic fracture through the AVI. The reader interested in the overall algorithmic implementation can directly go to Algorithm \ref{algAVI}.  In essence, we derive the proposed formulation from the discrete Hamilton's principle of stationary action with the fracture phase field incorporated. In addition, we adapt the reduced-space active set method to enforce the irreversibility constraint involved in the phase field problem. 

\subsection{Asynchronous discretization}
\begin{figure}[htbp]
\centerline{\includegraphics[width=0.7\textwidth]{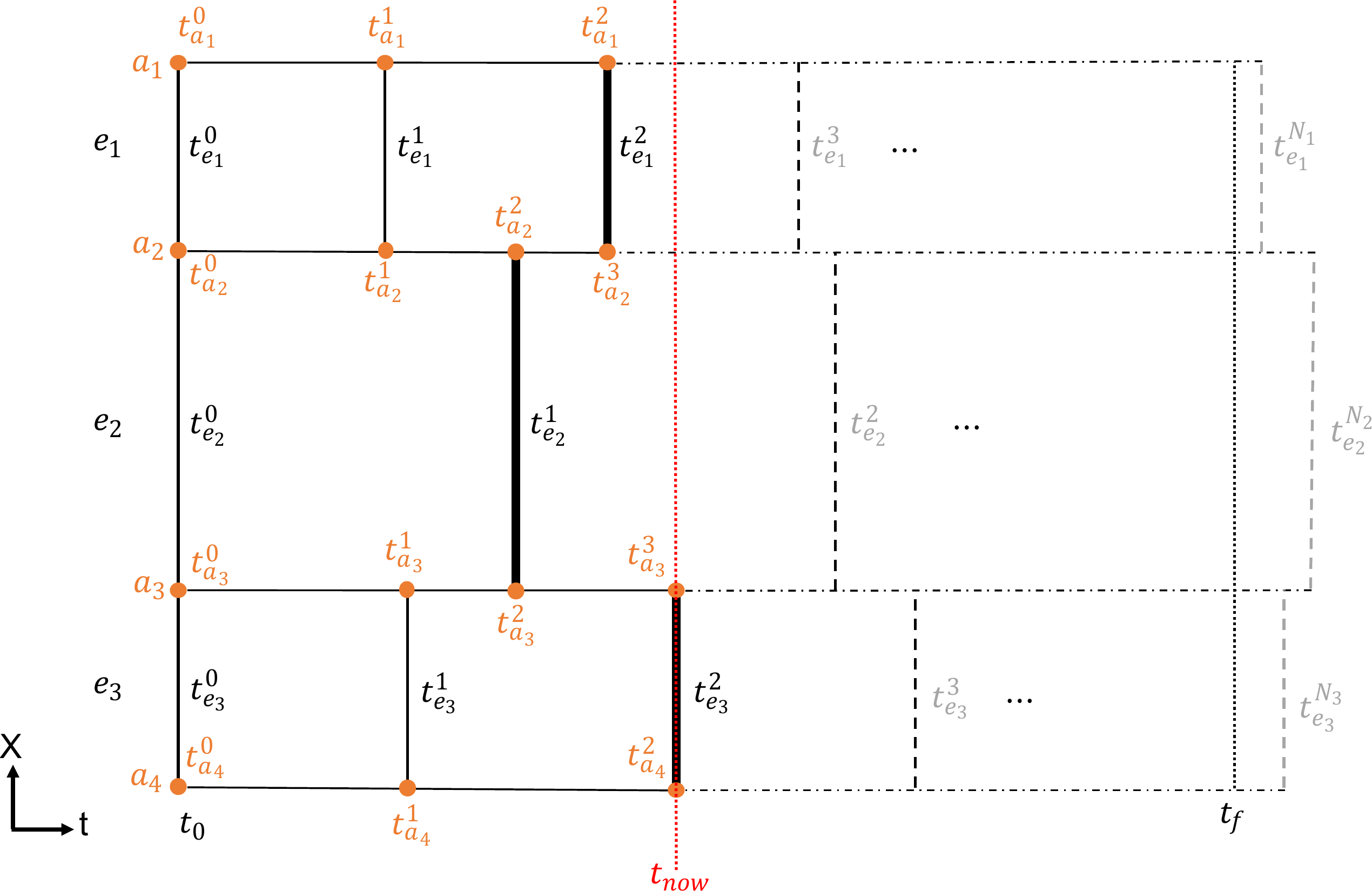}}
\caption{Asynchronous discretization of the time domain for a three-element mesh in the reference configuration. The entire update schedule follows the chronological order, i.e., $ \Theta = \left\{ t_{e_1}^0, t_{e_2}^0, t_{e_3}^0, t_{e_1}^1, t_{e_3}^1, t_{e_2}^1, t_{e_1}^2, t_{e_3}^2, t_{e_1}^3, t_{e_2}^2, t_{e_3}^3, \cdots \right\}$, where $ t_{now} = t_{e_3}^2$ is the current time and $e_3$ is the current active element.    
\label{fig:timedis}}
\end{figure}

This subsection presents the discretization of the time domain through an asynchronous strategy. Such asynchronous discretization allows each element to have an independent time step. As an example, 
Figure \ref{fig:timedis} shows the spacetime diagram of a three-element mesh with asynchronous time steps. 

Here, we assign to the element $ e \in \mathcal{T}_h $  the update schedule 
\begin{equation}\label{Eq:Theta_e}
	\begin{aligned}
		\Theta_e=\biggl\{t_0=t_e^0 < t_e^1 < t_e^2 \dots  t_e^{N_e-1} < t_f \le t_e^{N_e} \biggr\}.
	\end{aligned}
\end{equation}
At these instants, the displacements and velocities, and phase field values of all nodes $a\in \eta(e)$ are updated. 
In addition, we define the discrete elemental displacements $\mathbf{u}_e^j   \equiv\mathbf{u}_e(t_e^j)$ and discrete elemental phase fields $\mathbf{d}_e^j \equiv\mathbf{d}_e(t_e^j)$  at $t_e^j\in\Theta_e$, and the entire update schedule of the mesh is
\begin{equation}\label{Eq:Theta}
	\begin{aligned}
		\Theta=\bigcup_{e\in\mathcal{T}_h}\Theta_e.
	\end{aligned}
\end{equation}
For simplicity, we assume that there are no coincident instants  for any pair of elements except for the initial time, i.e., $\Theta_e\cap\Theta_{e'}=\{t_0\}$ if $e\neq e'$. The general case with coincident update instants can be handled without much difficulty and will not change the results as long as elements with coincident instants are far away enough from each other.

Similarly, we also gather the  schedule for node $a\in \eta $ as  
\begin{equation}\label{Eq:Theta_a}
	\begin{aligned}
		\Theta_a=\bigcup_{e\ni a} \Theta_e= \left\{t_0=t_a^0 < t_a^1 < \dots < t_a^{n_a-1}< t_a^{n_a} \right\}.
	\end{aligned}
\end{equation}
We additionally define  $\mathbf{u}_a^i \equiv  \mathbf{u}_a(t_a^i)$, $d_a^i\equiv d_a(t_a^i)$,
$t_a^i\in\Theta_a$, and the set of nodal displacements
\begin{equation}\label{Eq:Xi}
	\begin{aligned}
		\mathcal{U}_a=\left\{\mathbf{u}_a^i:i=0,1,\dots,n_a\right\},~a\in\eta,
	\end{aligned}
\end{equation}
and the set of nodal phase fields
\begin{equation}
    \mathcal{D}_a = \left\{d_a^i:i=0,1,\dots,n_a\right\},~a\in\eta.
\end{equation}

The triple $(\Theta,\mathcal{U}_a,\mathcal{D}_a)$ defines the discrete trajectory of the system. To solve for this triple, we write the discrete action sum as
\begin{equation}\label{Eq:action_disc}
	\begin{aligned}
		S_\mathrm{dis}(\Theta,\mathcal{U}_a,\mathcal{D}_a)=\sum_{e\in\mathcal{T}_h} \sum\limits_{j=0}^{N_e-1} L_e^j,
	\end{aligned}
\end{equation}
where $L_e^j \approx \int_{t_e^j}^{t_e^{j+1}} L_e\;\mathrm{d}t$.

This approximation can be realized by multiple schemes. In this paper,  we adopt one such that each node $a\in\eta$ follows a linear trajectory within the time interval $[t_a^i,t_a^{i+1}]$; consequently, the corresponding nodal velocities are constant in the said interval. Moreover, the potential energy and the crack energy terms are approximated with the rectangular rule using their values at $t_e^{j+1}$. Then the discrete Lagrangian is
\begin{equation}\label{eq:Lagrangian_disc}
		\int_{t_e^j}^{t_e^{j+1}}L_e \; \mathrm{d}t\approx L_e^j =\underset{a\in e}{\sum}\underset{t_a^i\in\left[t_e^j, t_e^{j+1}\right)}{\sum}\frac12  \mathbf{m}_{e,a}\bigl(t_a^{i+1}-t_a^i\bigr)\left \| \frac{\mathbf{u}_a^{i+1} - \mathbf{u}_a^i }{ t_a^{i+1} - t_a^i } \right \| ^2-\bigl(t_e^{j+1}-t_e^j\bigr)\left( V_e\bigl(\mathbf{u}_e^{j+1},\mathbf{d}_e^{j+1}\bigr) + \Gamma_{e}\bigl(\mathbf{d}_e^{j+1}\bigr)  \right),
\end{equation}
where $ \mathbf{m}_{e,a} $ is the mass matrix entry of node $a$ contributed by element $e$ and $\mathbf{d}_e^{j+1} \equiv \mathbf{d}_{e}(t_e^{j+1})$ is the elemental phase field vector. 
Finally, the discrete action sum  \eqref{Eq:action_disc} takes the following form
\begin{equation}\label{Eq:disc_action}
       S_\mathrm{dis} = \underset{a\in\eta}{\sum} \sum\limits_{i=0}^{n_a-1} \frac12 \mathbf{M}_a \bigl(t_a^{i+1}-t_a^i\bigr)\left \| \frac{\mathbf{u}_a^{i+1} - \mathbf{u}_a^i }{ t_a^{i+1} - t_a^i } \right \| ^2
        - \underset{e\in\mathcal{T}_h}{\sum} \sum\limits_{j=0}^{N_e-1}    \bigl(t_e^{j+1}-t_e^j\bigr) \left( V_e\bigl(\mathbf{u}_e^{j+1},\mathbf{d}_e^{j+1}\bigr) + \Gamma_{e}\bigl(\mathbf{d}_e^{j+1}\bigr)  \right),
\end{equation}
where $ \mathbf{M}_a  = \sum_{e\in\eta^{-1}(a)}\mathbf{m}_{e,a} $.

\subsection{Discrete variational principle}
In this subsection, we derive the formulation of the AVI for the phase field approach to dynamic fracture using  the discrete Hamilton's principle \cite{marsden_2001}. Taking the partial derivative of the discrete action sum  \eqref{Eq:disc_action} with respect to $\mathbf{u}_a^{i}$ follows
\begin{equation*}
    0 = \frac{\partial S_{dis}}{\partial \mathbf{u}_a^{i}} = \frac{\partial}{\partial \mathbf{u}_a^{i}} \left( \frac12 \mathbf{M}_a \bigl(t_a^{i}-t_a^{i-1}\bigr)\left \| \frac{\mathbf{u}_a^{i} - \mathbf{u}_a^{i-1} }{ t_a^{i} - t_a^{i-1} } \right \| ^2 + \frac12 \mathbf{M}_a \bigl(t_a^{i+1}-t_a^i\bigr)\left \| \frac{\mathbf{u}_a^{i+1} - \mathbf{u}_a^i }{ t_a^{i+1} - t_a^i } \right \| ^2 - \big(t_e^j - t_e^{j-1}\big) V_e\left(\mathbf{u}_e^{j},\mathbf{d}_e^{j}\right) \right),
\end{equation*}
where $a \in \eta(e)$ such at $t_a^i = t_e^j $, which yields the  discrete  Euler-Lagrange equations
\begin{equation}\label{eq:dis Euler}
    \mathbf{p}_a^{i+1/2} - \mathbf{p}_a^{i-1/2} = I_{e,a}^i = - \left(t_e^j - t_e^{j-1}\right) \frac{\partial V_e\left(\mathbf{u}_e^{j},\mathbf{d}_e^{j}\right)}{\partial\mathbf{u}_a^{i}},
\end{equation}
where $I_{e,a}^i$ may be regarded as the impulse component of node $a \in e $ at the time $t_a^i = t_e^j$,  and the discrete linear momentum is defined as
\begin{equation} \label{eq:p_a momenta}
        \mathbf{p}_a^{i-1/2}  =  \mathbf{M}_a \frac{\mathbf{u}_a^{i} - \mathbf{u}_a^{i-1} }{ t_a^{i} - t_a^{i-1} } = \mathbf{M}_a \mathbf{v}_a^{i-1/2}.
\end{equation}
Similarly, we take the partial derivative of \eqref{Eq:disc_action} with respect to $\mathbf{d}_e^{j}$ as follows 
\begin{equation}\label{stationarity}
     0 = \frac{\partial S_{dis}}{\partial \mathbf{d}_e^{j}} = \frac{\partial}{\partial \mathbf{d}_e^{j}} \left[ V_e\big(\mathbf{u}_e^{j},\mathbf{d}_{e}^{j}\big) + \Gamma_{e}\big(\mathbf{d}_{e}^{j}\big)  \right],
\end{equation} 
for element $e$ at time $ t_e^{j} $,  then the phase field of element $e$ are updated by
\begin{equation}\label{Eq:update_d}
        \mathbf{d}_e^{j}=\underset{\mathbf{d}_e^{\star}\leq \mathbf{d}_{e}\leq 1}{\argmin} \left\{ V_e\left(\mathbf{u}_e^{j},\mathbf{d}_{e}\right) + \Gamma_{e}\left(\mathbf{d}_{e}\right) \right\},
\end{equation}
where $\mathbf{d}_e^{\star}$ represents the phase field value of $\eta(e)$ at their most recent time of update; namely, for node $a\in \eta(e)$, at  step $i$,  $\mathbf{d}_e^{\star}$ contains the phase field value $d_a^{i-1}$. Note that normally in the same element $e$, $\mathbf{d}_e^{\star}$ may contain phase field values at different times. Here the stationarity condition \eqref{stationarity} becomes a minimization in \eqref{Eq:update_d} since $\Gamma_e$ is elliptic. The irrevesibility constraint in \eqref{Eq:update_d} may be enforced in many ways, for which we have chosen the reduced-space active set strategy, to be discussed in Section \ref{sec:activeset}.  

Now, we consider an element $e \in \mathcal{T}_h $ with $\mathbf{u}_{e}^{j-1}$ and $\mathbf{p}_a^{i-1/2}$, $a\in\eta(e)$, known at time $t_e^{j-1}$, also known $\mathbf{d}_e^{\star}$, the provisional solution procedure is thus:
\begin{itemize}
\item  For all $a\in\eta(e)$, solve $\mathbf{u}_a^{i}$ from \eqref{eq:p_a momenta}:~$ \mathbf{u}_a^{i} = \mathbf{u}_a^{i-1} + (t_a^{i}-t_a^{i-1}) \mathbf{M}_a^{-1} \mathbf{p}_a^{i-1/2} $.
\item Solve $\mathbf{d}_e^{j}$ from \eqref{Eq:update_d}.
\item  For all $a\in\eta(e)$, solve $\mathbf{p}_{a}^{i+1/2}$  from \eqref{eq:dis Euler}:~$\mathbf{p}_{a}^{i+1/2} = \mathbf{p}_a^{i-1/2} - \big(t_e^{j} - t_e^{j-1}\big) {\partial V_e\big(\mathbf{u}_e^{j}, \mathbf{d}_e^{j}\big)}/{\partial\mathbf{u}_a^{i}} $.
\end{itemize}

\subsection{Reformulation for solving the phase field with element patches}
\begin{figure}[htbp]
\centerline{\includegraphics[width=0.5\textwidth]{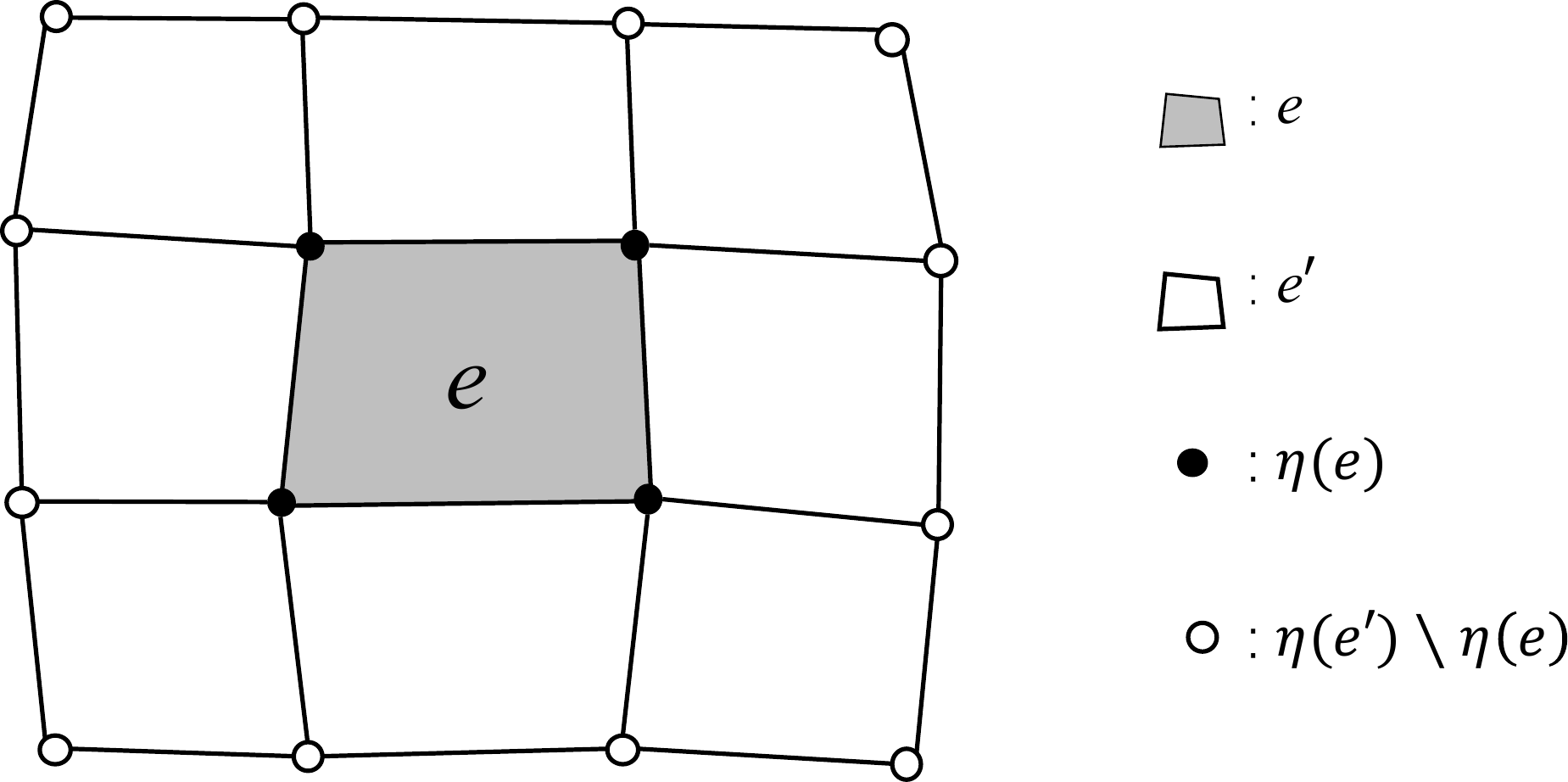}}
\caption{Diagram of a patch $\mathcal{T}_e$ that consists of the element $e$ (gray element) and its  adjacent elements $e'\in\mathcal{T}_e$ (white elements). At the beginning of an iteration for element $e$ at time $t_e^j$ for solving the phase field $\mathbf{d}_e^{j}$, the displacements $\mathbf{u}_e^j$ are known. Correspondingly, both the displacements $\mathbf{u}_{e'}^*$ and the phase fields $\mathbf{d}_{e'}^*$ of nodes in $\eta(e')\setminus\eta(e)$  (hollow nodes) assume their most recent values for the iteration.}
\label{fig:patch}
\end{figure}

The solution procedure mentioned above is variational; however, the results obtained with \eqref{Eq:update_d} show an unreasonable crack pattern (see Appendix \ref{app1}), hence we reformulate the constrained optimization problem \eqref{Eq:update_d} as follows. Essentially we want to minimize $V+\Gamma$ with the newly obtained $\mathbf{u}_e^{j}$ (same as before) while the field values of all nodes not belonging to $e$ frozen to their most recent values.
To this end, we define the patch for element $e$
\begin{equation}
    \mathcal{T}_e = \left\{ e'\in \mathcal{T}_h : \eta(e)\cap\eta(e')\neq\emptyset \right\},
\end{equation}
as shown in Figure \ref{fig:patch}. In this way,  Eq.~\eqref{Eq:update_d} is modified to take into account the contributions of its neighboring elements
\begin{equation}\label{Eq:patch_d}
        \mathbf{d}_e^{j}=\underset{\mathbf{d}_{e}^{\star }\leq \mathbf{d}_{e}\leq 1}\argmin 
        \sum_{\mathcal{T}_e} \left[ V_{e}\left( \big\{\mathbf{u}_{e'}^{*},\mathbf{u}_{e}^{j} \big\},\big\{\mathbf{d}_{e'}^{*},\mathbf{d}_{e} \big\} \right) + \Gamma_{e}\left(\left\{\mathbf{d}_{e'}^{*},\mathbf{d}_{e} \right\}\right) \right],
\end{equation}
where the superscript $*$ represents the nodal values of $\eta(e')\setminus\eta(e)$  (hollow nodes in Figure \ref{fig:patch}) at their most recent time of update.

Based on the spatial discretization, the minimization problem \eqref{Eq:patch_d} leads to the phase field residual of the element 
\begin{equation}\label{eq:r_patch}
    \mathbf{r}_e := \int_{\mathcal{T}_e} \left[ g'(d)\psi_+(\boldsymbol{\varepsilon}) N_a  + \frac{g_c}{4c_w}\left(\frac{w'(d)N_a}{\ell} + 2\ell \nabla d \cdot \nabla N_a \right)  \right] \mathrm{d}\Omega,\quad a\in\eta(e),
\end{equation}
where $g'=\mathrm{d}g/\mathrm{d}d$ and the tangent stiffness matrix of the element 
\begin{equation}\label{eq:w_patch}
    \mathbf{k}_e := \int_{\mathcal{T}_e} \left[  g''(d)\psi_+(\boldsymbol{\varepsilon}) N_a N_{a'} + \frac{g_c}{4c_w} \left(\frac{w''(d)N_a N_{a'} }{\ell} + 2\ell \nabla N_a \cdot \nabla N_{a'} \right) \right] \mathrm{d}\Omega,\quad a,a'\in\eta(e).
\end{equation}
For the detailed derivation of \eqref{eq:r_patch} and \eqref{eq:w_patch}, the reader is referred to Shen et al.~\cite{shen2018implementation}

\subsection{Reduced-space active set method for irresversibility constraint} \label{sec:activeset}

There are several approaches to impose the inequality constraints of the phase field when solving \eqref{Eq:patch_d}, such as the local history variable method\cite{MIEHE20102765}, the penalty method\cite{Gerasi2019Penalization}, and the augmented Lagrangian method\cite{GEELEN2019680}. 

In this work, we employ the reduced-space active set strategy\cite{WU2018Robust} to ensure the phase field bounds $d\in[0,1]$ and the irreversibility condition $\dot d>0$. In the discrete setting, the phase field needs to satisfy the condition
\begin{equation}
    0 \le d_a^{i-1} \le d_a^i \le 1, \quad \forall i = 1,2,\cdots, n_a \text{ and } a\in\eta(e).
\end{equation}
Note that we solve this inequality-constrained optimization problem efficiently for only one element instead of for the entire domain. Then the solutions are determined by a mixed complementarity problem\cite{Farrell2017Linear}  
\begin{equation}\label{eq:mix_condition}
        \left\{
    \begin{aligned}
     & d_a^{i-1} = d_a^i,  & r_a \ge 0, \\ 
     & d_a^{i-1} \le d_a^i \le 1, & r_a = 0, \\
     & d_a^i = 1, & r_a \le 0,
     \end{aligned}
    \right.
\end{equation}
where $d_a^i$ is the nodal phase field  at time $t_a^i = t_e^j$, and $r_a$ is the phase field residual corresponding to node $a$. For each iteration, with  $\mathbf{d}_e^{\star}$, $\mathbf{u}_e^{j}$, $\mathbf{d}_{e'}^{*}$, and $\mathbf{u}_{e'}^{*} $ of the patch at time $t_e^j$ known, the phase field value of element $e$ can be updated through Algorithm \ref{Algo: active set}. 

\begin{algorithm}
\caption{Reduced-space active set method for solving the phase field}\label{Algo: active set}
\hspace*{0.02in} {\bf Input:} $ \mathbf{d}_e^{\star}, \mathbf{u}_e^{j}, \mathbf{d}_{e'}^{*}, \mathbf{u}_{e'}^{*} $, $\forall e'\in\mathcal{T}_e$   \\
\hspace*{0.02in} {\bf Output:} $ \mathbf{d}_e^{j}$
\begin{algorithmic}[1]
  \State Compute $\mathbf{r}_e $ and define flag $\leftarrow$ true
   \If {$\|\mathbf{r}_e\|$ < tol }
   \State flag $\leftarrow$ false
   \Else
    \State Define $  \mathcal{A} := \left\{ a\in \eta(e):  \left(d_a^i = d_a^{i-1}, r_a > 0\right) \text{ or } \left(d_a^i  = 1, r_a<0\right) \right\} \text{and } \mathcal{I} = \eta(e)\setminus\mathcal{A}$
   \EndIf
   \While {flag == true}
   \State Compute $ \mathbf{k}_e $ and let  $\mathbf{d}_{\mathcal{A}} \leftarrow \mathbf{d}_{\mathcal{A}}^{j-1}$,  $\mathbf{d}^{0}_{\mathcal{I}} \leftarrow \mathbf{d}_{\mathcal{I}}^{j-1} $
    \While {$\|\mathbf{r}_{\mathcal{I}}\|$ > tol }  
    \State Solve  $\mathbf{d}_{\mathcal{I}}^{n+1} \leftarrow \mathbf{d}_{\mathcal{I}}^{n} - \mathbf{k}_{\mathcal{I}}^{-1} \mathbf{r}_{\mathcal{I}}$ and compute $\mathbf{r}_{\mathcal{I}}$, then let $n \leftarrow n+1$
   \EndWhile
   \For {$a\in \mathcal{I} $}
     \If {$d_a^{n+1} >1 $ } 
     \State Set $d_a^{n+1} \leftarrow 1$ and $\mathcal{A} \leftarrow \mathcal{A} \cup\{a\}$, $\mathcal{I} \leftarrow \mathcal{I}\setminus\{a\}$
     \ElsIf{$d_a^{n+1} < d_a^{n}$}
     \State Set $d_a^{n+1} \leftarrow d_a^{n}$ and $\mathcal{A} \leftarrow \mathcal{A} \cup\{a\}$, $\mathcal{I} \leftarrow \mathcal{I}\setminus\{a\}$
    \EndIf
   \EndFor
    \State Compute $\mathbf{r}_{\mathcal{A}}$ from \eqref{eq:r_patch}
    \If {$\forall a \in \mathcal{A} $ satisfy  Eq.~\eqref{eq:mix_condition}}
    \State flag $\leftarrow$ false
    \Else 
    \State For each $ a \in \mathcal{A} $ not satisfying  Eq.~\eqref{eq:mix_condition}, $\mathcal{I} \leftarrow \mathcal{I} \cup\{ a\}$
    \EndIf
   \EndWhile
\end{algorithmic}
\end{algorithm}

\subsection{Algorithmic implementation}
This section focuses on the algorithmic implementation of AVI for the phase field fracture. The overall pseudo-code is provided in Algorithm \ref{algAVI}. The time step of each element is taken as a fraction of their critical time step and is computed by
\begin{equation}
    t_{crit} =  C_{CFL}\;\frac{2 }{\omega_e},
\end{equation}
where $C_{CFL}$ is taken as 0.6 and  $\omega_e$ is the maximum natural frequency of the element, which is the square root of the maximum eigenvalue of the generalized eigenvalue problem $ \mathbf{k}_e\mathbf{U}=\omega^2 \mathbf{m}_e \mathbf{U} $. The time step of each element allows certain adaptivity, although we keep $C_{CFL}$ constant in this work for simplicity. 

Due to the asynchrony of the algorithm, we employ the priority queue \cite{Knuth1998Art} to keep track of the causality. The priority queue assigns each element a priority according to their next update time where the element to be updated at a sooner time has a higher priority. In other words, the priority queue ensures that all elements in the queue are ordered according to their next time to be updated, and the top element in the queue is always the one whose next update time is the closest to the current time in the future. 

The implementation details are shown in Algorithm \ref{algAVI}.  First, the first time steps of all elements in the mesh are computed and pushed into the priority queue to establish the initial queue. Within each iteration, the priority queue pops an element (calls the active element) and its next update time. The nodal displacements, phase fields, and momenta of the active element are updated accordingly. Subsequently, the next update time of this element is computed and if this time is less than $t_f$, it is pushed into the priority queue. The algorithm continues until the priority queue is empty.

\begin{algorithm}
\caption{Algorithm of AVI for the phase field to fracture}
\label{algAVI}
\hspace*{0.02in} {\bf Input:} $ \mathcal{T}_h,~\eta,~\Theta$, $Q = \emptyset $ and $\{ \mathbf{u}_a^0,~ \mathbf{d}_e^0,~ \mathbf{p}_a^{1/2}~ | a\in\eta\}$ \\
\hspace*{0.02in} {\bf Output:} $ \mathbf{u}_a^i, \mathbf{d}_e^j, \mathbf{p}_a^{i+1/2} $, where $i,j = 1,2,3,\cdots,$ corresponding to $\Theta$
\begin{algorithmic}[1]
   \State Initialization: $ \mathbf{u}_a \leftarrow \mathbf{u}_a^0,~\mathbf{d}_e \leftarrow \mathbf{d}_e^0,~\mathbf{p}_a \leftarrow \mathbf{p}_a^{1/2},~\tau_a \leftarrow t_0 $ for all $a \in \eta$
   \For {all $e\in \mathcal{T}_h$}
    \State $\tau_e \leftarrow t_0$
    \State Compute $t_e^1$ and push $(t_e^1,e)$ into priority queue $Q$
   \EndFor 
   \While {$Q$ is not empty}
   \State Extract next element: pop $(t,e)$ from $Q$
   \State Compute displacements $\mathbf{u}_e$ with \eqref{eq:p_a momenta} and update node's time $\tau_a \leftarrow t$ for all $a\in \eta(e)$ 
    \State Compute the phase fields $\mathrm{d}_e$ with \eqref{Eq:patch_d}  following  Algorithm \ref{Algo: active set} 
    \If {$t<t_f$}
    \State Compute momentum $\mathbf{p}_a$ with \eqref{eq:dis Euler} for all $a\in \eta(e)$
    \State Update element's time: $\tau_e \leftarrow t$ 
    \State Compute $t_e^{next}$  and schedule $e$ for next iterate: push $(t_e^{next},e)$ into $Q$
    \EndIf
   \EndWhile
\end{algorithmic}
\end{algorithm}

\section{Numerical examples}\label{sec4}
In this section, we showcase three benchmark examples  to demonstrate the ability of the proposed formulation in capturing the key features of dynamic fracture. In addition, we compare the computational costs and the energy conservation behavior of our approach for AT1 and AT2 models. In all examples, we use the unstructured first-order quadrilateral finite elements  and refine the mesh along the potential crack path. 

A note on the post-processing is as follows. We sample the solution at a frequency of every 500,000 elemental iterations. For example, if the current time $t_{now}=t_{e_3}^2$ in Figure \ref{fig:timedis} happens to be a sampling time, then post-processing results are obtained using the most recent nodal displacement, velocity, and phase field values prior to $t_{now}$, i.e., values at nodal times $t_{a_1}^2, t_{a_2}^3, t_{a_3}^3, t_{a_4}^2$ for the nodes shown.  For example, the crack patterns to be plotted are obtained using the most recent phase field values prior to the sampling times.  More accurate results can be obtained by  interpolation using the values before and after the sampling time, which is not undertaken in this work for simplicity.

\begin{table}[tbp]%
\centering
\caption{Material properties for the numerical examples in Section~\ref{sec4}. \label{mater_pro}}%
\begin{tabular*}{500pt}{@{\extracolsep\fill}lcccc@{\extracolsep\fill}}
\toprule
 \textbf{Parameter}   & \textbf{Symbol} & \textbf{Section~\ref{subsec_branch} }  & \textbf{Section~\ref{sub_CTest} }  & \textbf{Section~\ref{sub_kalthf}}  \\
\midrule
Material  & -  &   Silica glass  & Soda-lime glass  & Maraging steel 18Ni(300) \\
Young's modulus (GPa) & $E$    & 32    & 72     & 190  \\
 Poisson's ratio       & $\nu$  & 0.2   & 0.22   & 0.3  \\
Density (kg/m$^3$)    & $\rho$ & 2450  & 2440   & 8000 \\
Critical energy release rate (J/m$^2$) & $g_c$ & 3  & 3.8    & 2.213$\times$10$^4$\\
 Rayleigh speed (m/s)      & $v_R$ & 2119  & 3172 & 2803 \\
\bottomrule
\end{tabular*}
\end{table}

\subsection{Boundary tension test}\label{subsec_branch}
\begin{figure}[htbp]
\centerline{\includegraphics[width=0.35\textwidth]{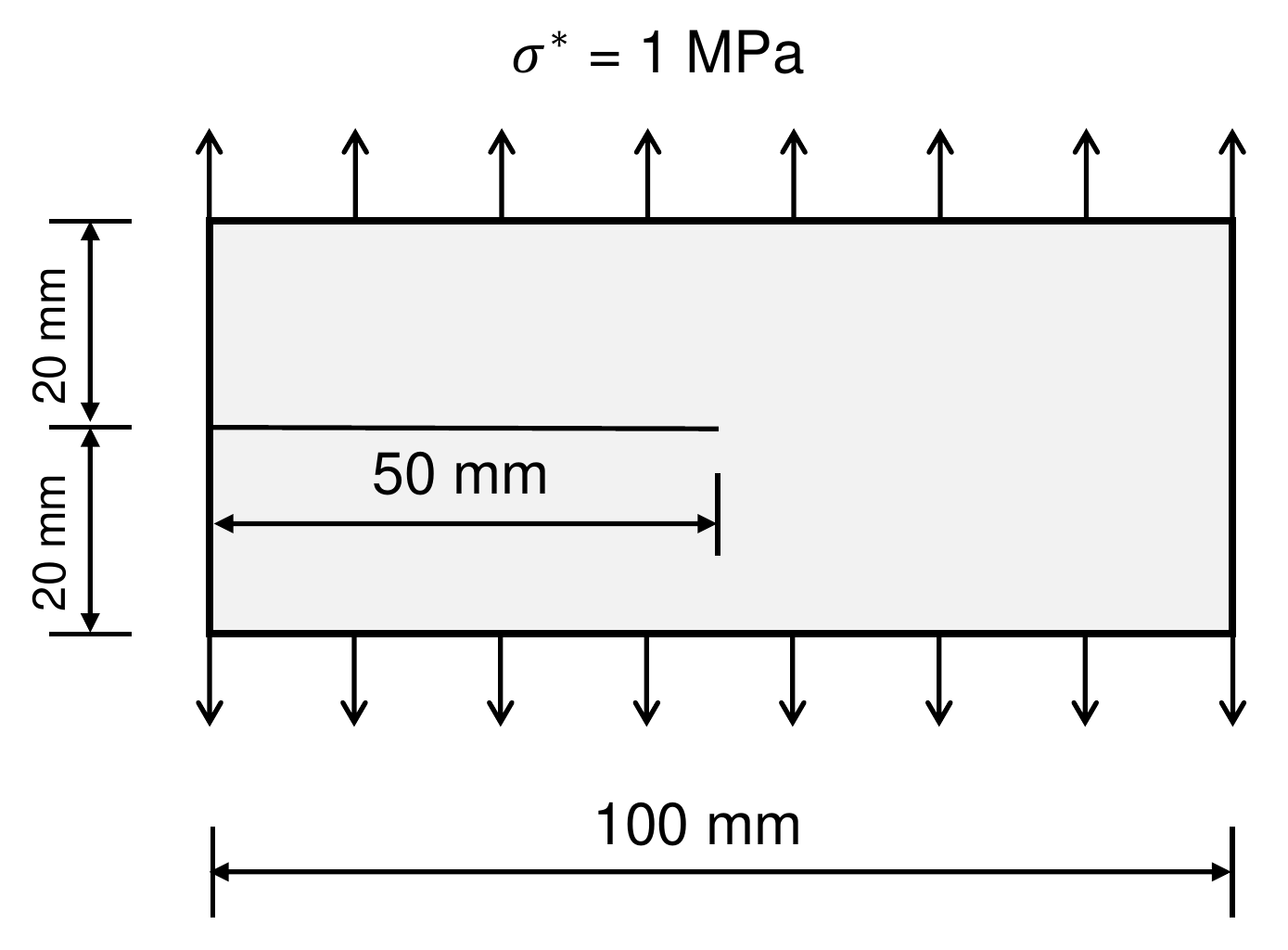}}
\caption{Geometry and boundary conditions of the boundary tension test, in which a pre-notched plate is under tension.  }\label{fig:geo_traction}
\end{figure}

In this section, a pre-notched rectangular plate loaded dynamically in tension is modeled. The geometry and boundary conditions are shown in Figure \ref{fig:geo_traction}. A constant traction $\sigma^{*} = \SI{1}{MPa}$ is applied on the top and bottom edges throughout the simulation and the remaining boundary is traction free. This benchmark problem has been widely studied, for example by Song et al.~\cite{Song2008_comparative} using the extended finite element method, by Nguyen \cite{Nguyen2014Discontinuous} with the cohesive zone method, and by Borden et al.~\cite{BORDEN201277} with a synchronous phase field approach to fracture, as well as in experimental studies \cite{Ramulu1985Mechanics,Sharon1996Micro}. As described in Ref.~\citenum{Ramulu1985Mechanics}, a crack  emerges at the notch tip and starts propagating to the right in a stable way. Over a certain distance, the  main crack branches into two symmetrical sub-cracks and continue growing until it reaches the right surface.

The material used in this test is silica glass and its properties are listed in Table \ref{mater_pro}. A plain strain state with a unit thickness is assumed. The length scale parameter takes $ \ell = 5\times10^{-4}$~m, which is small enough with respect to the specimen dimensions. Two different mesh levels are used: Mesh 1 with $ h= 2.5\times 10^{-4} $~m $=\ell/2$, and Mesh 2 with $ h= 1.25\times 10^{-4} $~m $=\ell/4$.

\begin{figure}[tbp]
    \centering
    \includegraphics[width=0.9\textwidth]{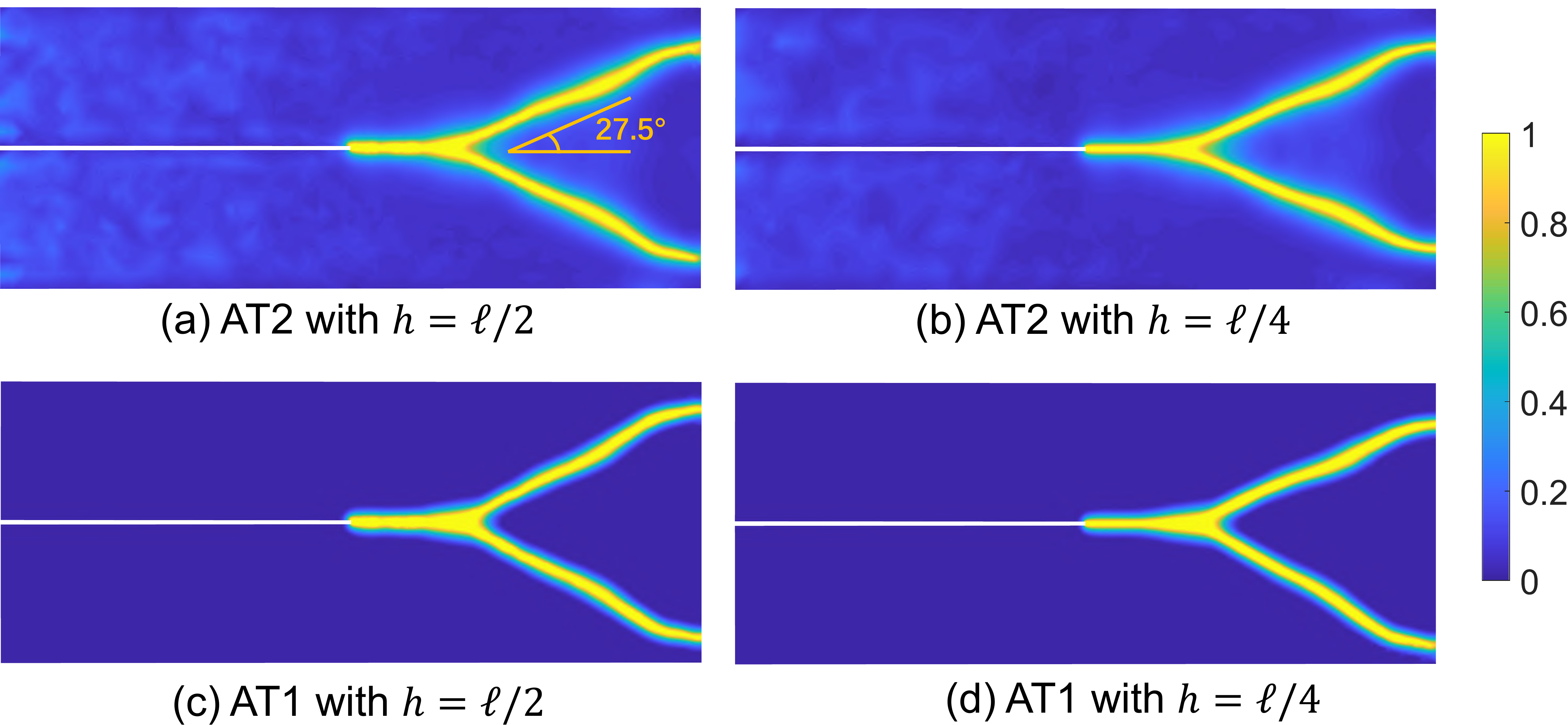}
    \caption{Phase field results for the boundary tension test of AT1/2 models.}
    \label{fig:traction_pf}
\end{figure}

The final phase field results are shown in Figure \ref{fig:traction_pf}. As seen, there is no significant difference of the crack pattern between the AT2 and AT1 models. The crack branches at between 34 and 36 $\mu$s and reaches the right boundary at $t \approx 80 \mu$s. 
The upper crack branching angle is around $27.5 ^\circ$, which  agrees well with the results in Refs.~\citenum{Schluter2014Phase, Vinh2018modeling}. In addition, the bifurcation angle of the lower branch is slightly different from that of the upper one, which may be caused by the non-symmetric discretization of the mesh. This non-perfect symmetry was also observed by Ren et al.~\cite{REN2019explict}.

\begin{figure}[tbp]
    \centering
    \includegraphics[width=0.95\textwidth]{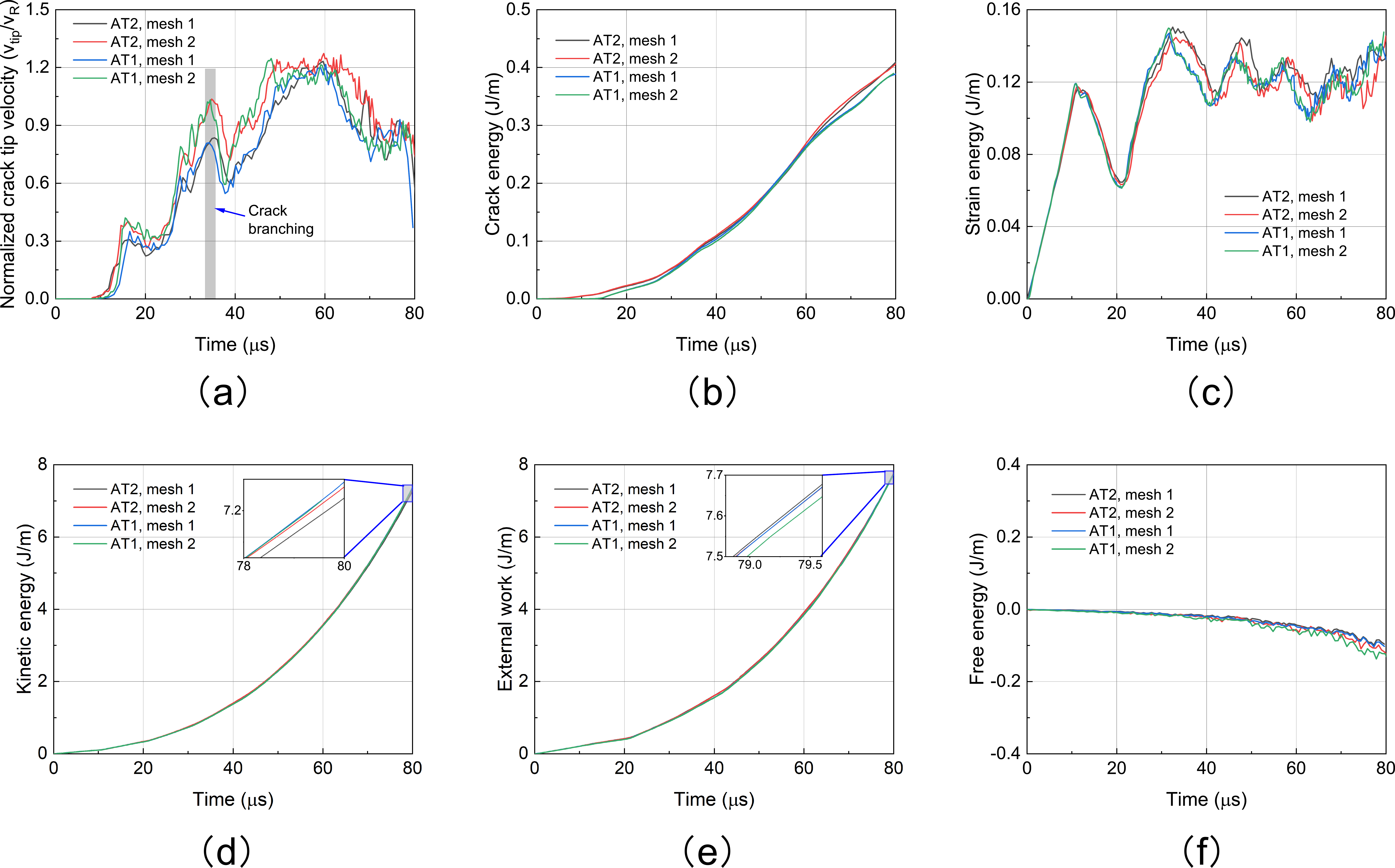}
    \caption{Results of the boundary tension test. Evolution of (a) normalized total crack tip velocity; (b) crack energy $\Gamma(d)$; (c) strain energy; (d) kinetic energy $T(\boldsymbol{\dot u})$; (e) external work; (f)  free energy $T(\boldsymbol{\dot u}) + V(\boldsymbol{u}, d) + \Gamma(d)$. Here for (e) as well as for (b)(c)(d) the results are obtained by sampling. It can be seen that the free energy is 1.32\% of the external work at the end,  indicating the conservation of energy.}
    \label{test1_vel_ene}
\end{figure}

Figure \ref{test1_vel_ene}(a) shows the evolution of the total crack tip velocity calculated  by 
\begin{equation}\label{eq:speed}
  v_{tip} = \frac{1}{g_c} \frac{\mathrm{d}\Gamma}{\mathrm{d} t},
\end{equation}
and normalized by the Rayleigh wave speed. In particular, the derivative is obtained by comparing the values of $\Gamma$ at consecutive sampling times. It is observed that at the beginning, a single main crack propagates to the right with a speed less than 60\% of the Rayleigh wave speed. Then, the main crack branches into two sub-cracks and in this respect the total crack tip velocity of both branches is plotted, which is still less than 60\% of twice the Rayleigh wave speed. Therefore, whether before or after the branching emerges, the velocity is in a reasonable range. Moreover, the overall propagation speed during the evolution is in good agreement with the results reported in a numerical study \cite{BOURDIN2000797} and an experimental study \cite{Ramulu1985Mechanics}.  

Figure \ref{test1_vel_ene} (b)-(d) present the evolution of the crack surface energy, strain energy, and kinetic energy, respectively.  The crack surface energy  monotonically increases as expected due to the unilaterality of the phase field. In addition, the strain energy evidently shows the  periodic oscillation at the beginning and this trend gradually weakens with crack evolution, because the stress wave is reflected at the boundaries and cracks, and interacts with itself.

Figure \ref{test1_vel_ene}(e) shows the evolution of the external work, which is calculated from the second term on the right hand side of Eq.~\eqref{eq:regular potential}. It is clear that the kinetic energy accounts for most of the  energy converted from external work. 

Figure \ref{test1_vel_ene}(f) shows the evolutions of the free energy $T(\boldsymbol{\dot u}) + V(\boldsymbol{u}, d) + \Gamma(d)$. The free energy is negative and its magnitude is only 1.32\% of the external work at the end. This small negative numbers demonstrate that the method possesses remarkable energy conservation property and is energetically stable.

\begin{figure}[htbp]
    \centering
    \includegraphics[width=0.8\textwidth]{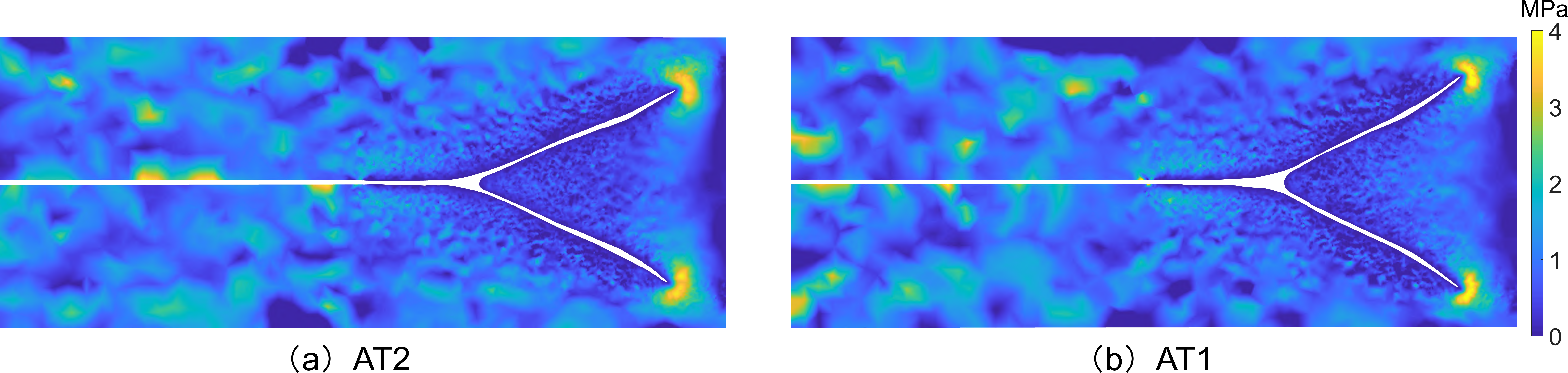}
    \caption{Boundary tension test: maximum principal stress with Mesh 1 at $t = 70$ $\mu$s, where $d >0.9$ are removed. }
    \label{fig:MaxStr_traction}
\end{figure}

Figure \ref{fig:MaxStr_traction} shows the maximum principal stress with Mesh 1 at $t = 70$ $\mu$s. Therein, stress concentration is clearly seen at the crack tips and the results are in good agreement with those in Ref.~\citenum{LIU2016Abaqus}.

Table \ref{tab:traction} collects some statistics of the computational cost for the example at hand. As a platform-independent indicator, the number of updates of each element throughout the simulation for each case is counted. The second, third, and fourth columns represent the  minimum, maximum, and median  numbers of updates among the elements, respectively. The fifth column is the total numbers of updates of all elements of AVI. The sixth column is the total numbers of updates of synchronous integration, where the data is estimated by assuming the global critical time step is used for the same time interval $[t_0, t_f]$.  As shown, the total numbers of AVI updates is approximately 31\%  of those of synchronous integration. Considering that it is even more costly to implicitly solve for the phase field with a synchronous method per time step, the data in  Table \ref{tab:traction} indicates that the proposed scheme effectively reduces the computational cost compared with synchronous methods.  

\begin{table}[htbp]%
\centering 
\caption{Numbers of elemental updates for the boundary tension test during the entire simulation. }
\label{tab:traction}
\begin{tabular*}{500pt}{@{\extracolsep\fill}lccccc@{\extracolsep\fill}}
\toprule
\textbf{Mesh \& Model} & \textbf{Min.}  & \textbf{Max.}  & \textbf{Median}  & \textbf{AVI total\tnote{$^\dagger$}}   & \textbf{Synchronous integration (estimated)\tnote{$^\ddagger$}} \\
\midrule
Mesh 1, AT2 & 178  & 5,983   & 2,106 & 27,400,002  & 86,382,554   \\
Mesh 1, AT1 & 178  & 5,983   & 2,106 & 27,400,002  & 86,382,554   \\
Mesh 2, AT2 & 160  & 10,150  & 3,440 & 62,285,189  & 201,051,200  \\
Mesh 2, AT1 & 160  & 10,150  & 3,440 & 62,285,189  & 201,051,200  \\
\bottomrule
\end{tabular*}
\begin{tablenotes}
\item[$^\dagger$]Total numbers of the elements involved in the update of the mechanical field and phase field.
\item[$^\ddagger$] This column of data is estimated by assuming the global critical time step is used throughout the computation for the same desired time interval $[t_0,t_f]$.
\end{tablenotes}
\end{table}

\subsection{Compact tension test}\label{sub_CTest}

\begin{figure}[htbp]
\centerline{\includegraphics[width=0.4\textwidth]{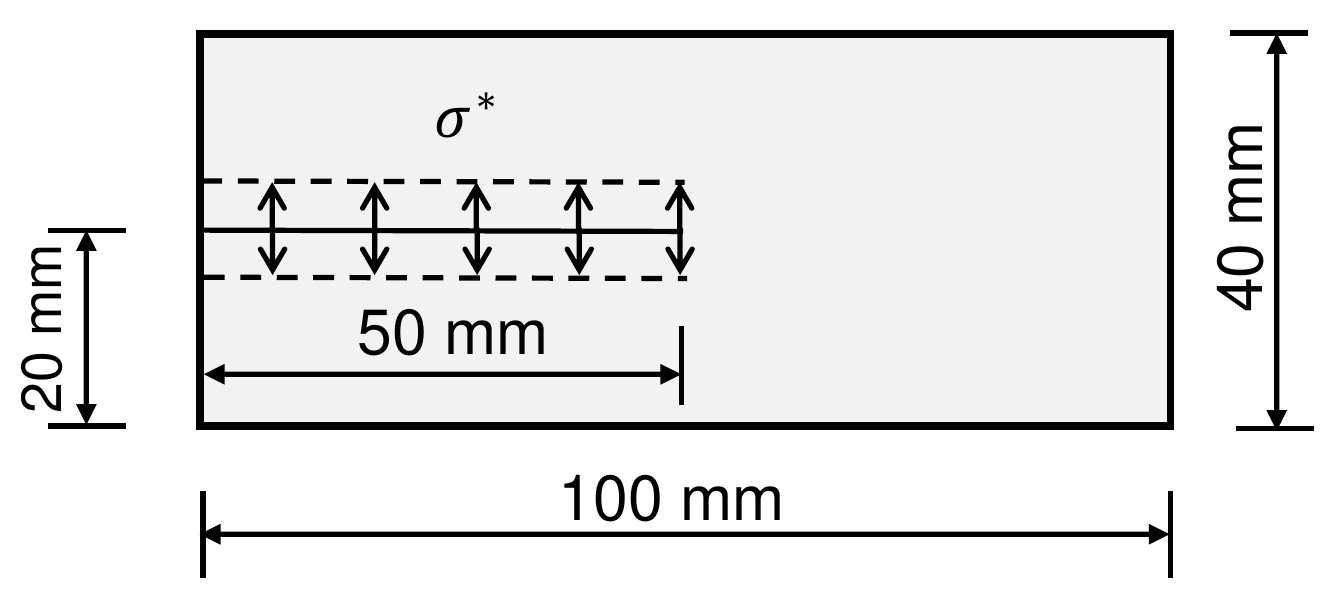}}
\caption{Geometry and boundary conditions for the dynamic CT test. \label{fig:geo_CTest}}
\end{figure}

In this section, we investigate a series of dynamic loads applied on pre-crack surfaces as the compact tension (CT) test. The geometry and boundary conditions are shown in Figure \ref{fig:geo_CTest}. Three different constant normal tractions $\sigma^* = \{0.5, 3, 6\}$~MPa are applied on the pre-crack surfaces. This benchmark problem has been studied by  Bobaru and Zhang \cite{Bobaru2015Why} using peridynamics and Mandal et al.~\cite{MANDAL2020evaluation} with a synchronous phase field approach.

\begin{figure}[htbp]
    \centering
    \includegraphics[width=0.85\textwidth]{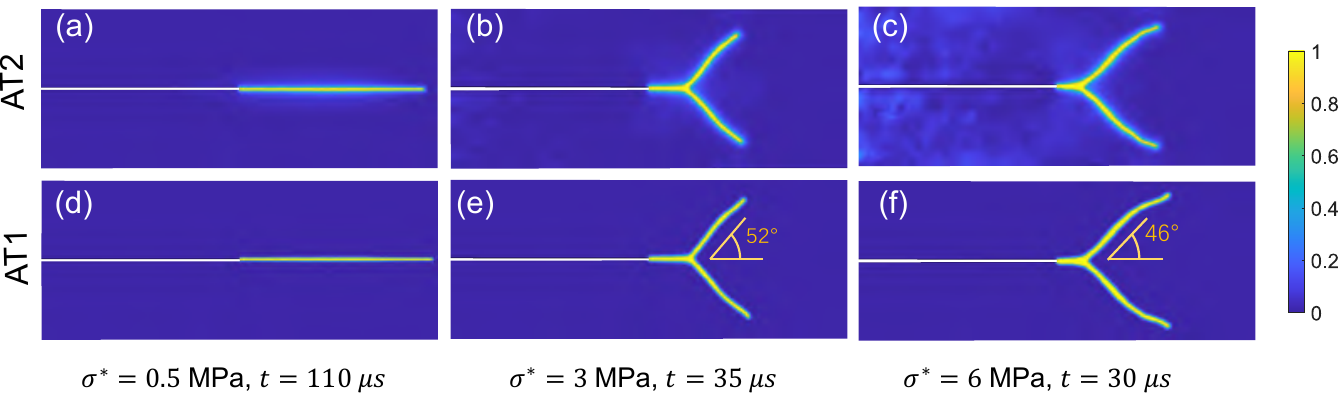}
    \caption{Phase field results for the CT test of AT1/2 models under different loads. }
    \label{fig:pf_CTset}
\end{figure}

The material is assumed to be soda-lime glass, whose properties are given in Table \ref{mater_pro}. Plain strain state is assumed. The length scale parameter $ \ell = 5\times10^{-4}$ m and the mesh size $ h= 2.5\times 10^{-4} $ m $=\ell/2$ are used for all cases.  Figure \ref{fig:pf_CTset} shows the phase field results for the CT test. For $ \sigma^* = 0.5$ MPa a straight crack without branching is obtained. For larger values of $\sigma^*$ crack branching is observed and the branching location moves to the left with the increase of $\sigma^*$. The crack branching happens at around  17.3 $\mu$s and 9.2 $\mu$s, and the branching angles are $52 ^\circ$ and $46 ^\circ$ for $\sigma^*=$ 3 MPa and 6 MPa, respectively. Also, there is no significant difference of the crack patterns between the AT2 and AT1 models for the same load. Moreover, the crack patterns, branching instants, and branching angles are all in good agreement with the results reported in Ref.~\citenum{MANDAL2020evaluation}.

\begin{figure}[tbp]
    \centering
    \includegraphics[width=0.95\textwidth]{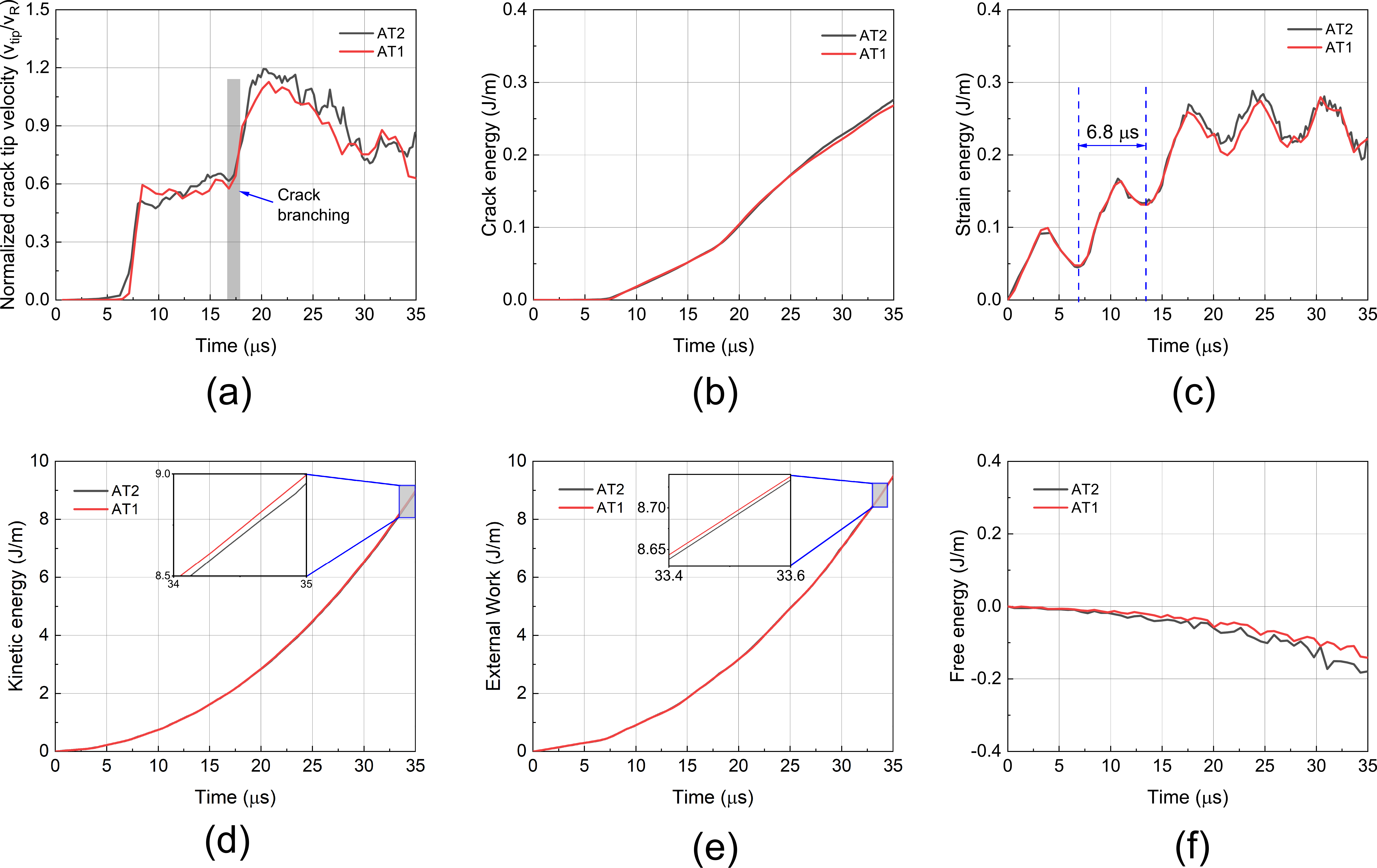}
    \caption{ Results of the  CT test for $\sigma^* = 3$ MPa. Evolution of: (a) normalized total crack propagation velocity; (b) crack energy; (c) strain energy; (d) kinetic energy; (e) external work; (f) free energy. The free energy is only  1.69\% of the external work, which indicates the conservation of energy. } 
    \label{fig:VSE_CTset}
\end{figure}

Figure \ref{fig:VSE_CTset}(a) illustrates the normalized total crack tip velocity of CT test for $\sigma^* = 3$ MPa. Like the case of Figure  \ref{test1_vel_ene}(a), a main crack propagates to the right with an increasing speed less than 60\% of the Rayleigh wave speed. Then, then main crack branches into two sub-cracks and the total speed is still less than 60\% of twice the Rayleigh wave speed. Figure \ref{fig:VSE_CTset}(b) shows the crack energy. Figure \ref{fig:VSE_CTset}(c) shows the evolution of the strain energy.  An interesting observation is that the curve presents a periodic oscillation with a period of approximately 6.8 $\mu$s, which can be explained as follows. During the process, the stress waves propagate from the crack to the top and bottom boundaries and then are reflected until they meet the crack again.  The time it takes the stress wave to travel a round trip can be estimated by $l/v_D = 6.89 $ $\mu$s, where $ l = 40$ mm is twice the half-width of the specimen and $v_D = 5800$ m/s is the dilatational wave speed of soda-lime glass. This process is repeated, and hence the periodicity.  Figure \ref{fig:VSE_CTset}(d-e) show the kinetic energy and external work, both of which monotonically increase.  Figure \ref{fig:VSE_CTset}(f) shows the free energy of the AT2 and AT1 model during the evolution.  As we can see, the magnitude of the free energy only accounts for 1.69\% of the external work, which indicates the conservation of energy. 

\begin{figure}[htbp]
    \centering
    \includegraphics[width=0.9\textwidth]{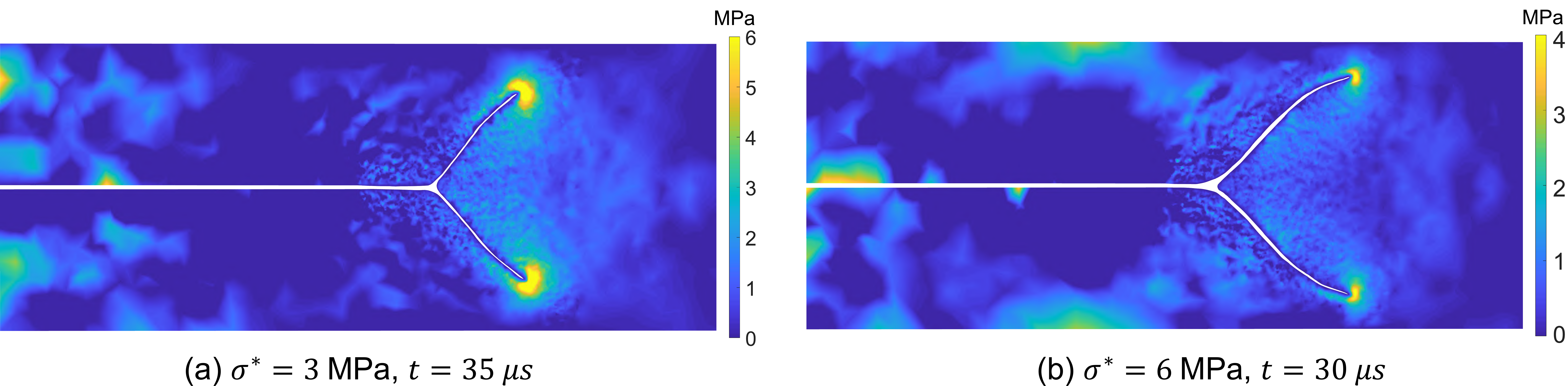}
    \caption{Maximum principal stress from the CT test for the AT1 model. }
    \label{fig:MaxStr_CTest}
\end{figure}

Figure \ref{fig:MaxStr_CTest} shows the maximum principal stress for $\sigma^* = 3$~MPa and 6~MPa, respectively. 

\subsection{The Kalthoff-Winkler test}\label{sub_kalthf}

\begin{figure}[htbp]
\centerline{\includegraphics[width=0.35\textwidth]{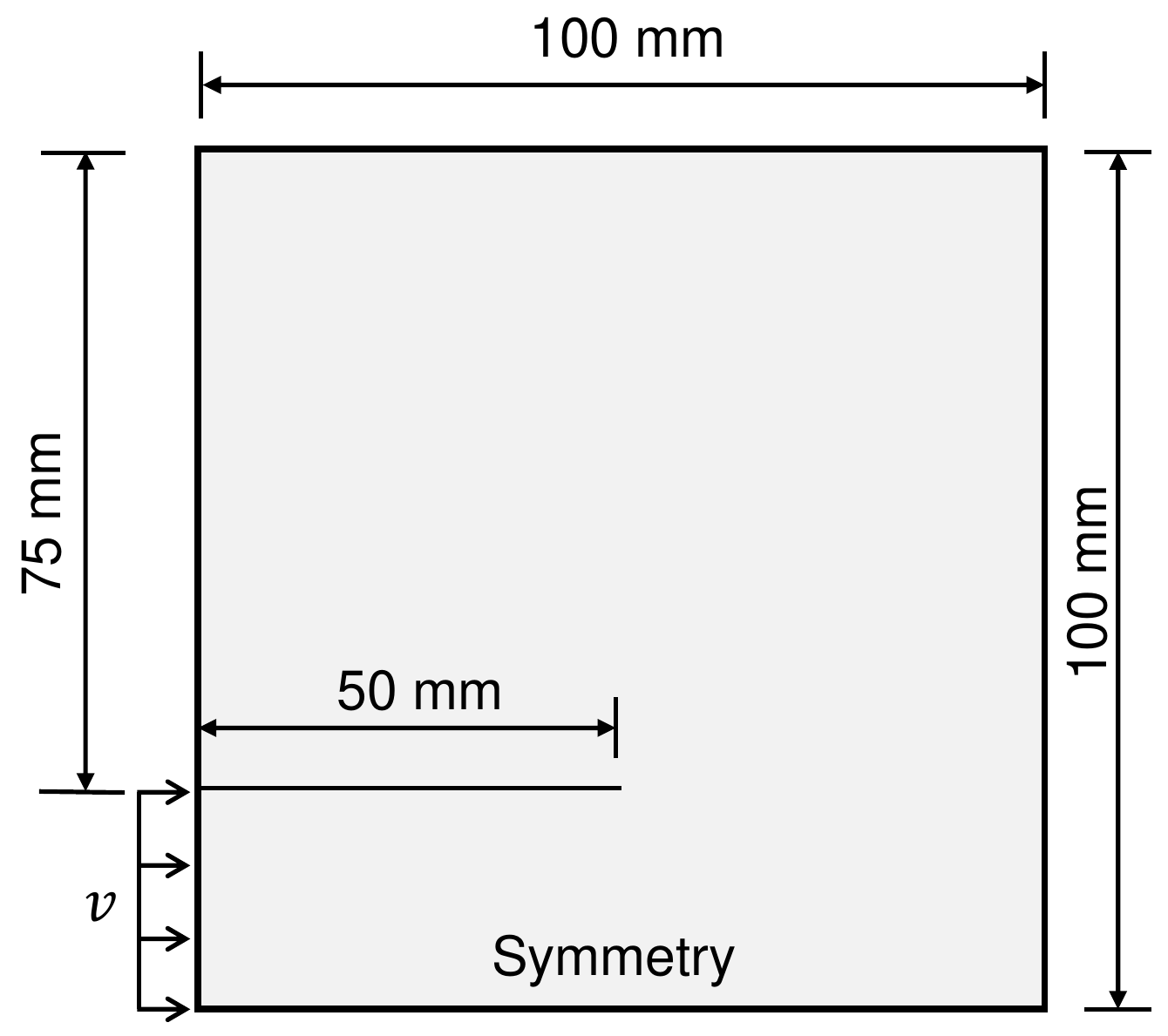}}
\caption{Geometry and boundary conditions of the Kalthoff test, where $v = 16.5 $ m/s.   \label{fig:geo_kalthoff}}
\end{figure}

This section studies the Kalthoff-Winkler experiment in which an edge-cracked plate is under impact velocity. Due to symmetry, only half of the plate is considered. The geometry and boundary conditions are shown in Figure \ref{fig:geo_kalthoff}. In the experiment \cite{kalthoff1988failure,kalthoff2000modes}, the brittle failure mode with a crack propagating at about $70 ^\circ$  was observed at a certain impact speed, and the relevant numerical results were reported by other researchers using the extended finite element method \cite{Song2008_comparative}, peridynamics \cite{WANG2019nonordinary}, and the gradient damage method \cite{Li2016Gradient}. 

The material is maraging steel 18Ni(300), whose properties are given in Table \ref{mater_pro}. A plane strain state is assumed. The length scale parameter $ \ell = 3.9\times10^{-4}$~m and two different meshes are used: Mesh 1 with size $ h= 1.95\times 10^{-4} $ m~$ =\ell/2$ and Mesh 2 with $ h= 9.75\times 10^{-5} $ m$=\ell/4$.

\begin{figure}[tbp]
    \centering
    \includegraphics[width=0.65\textwidth]{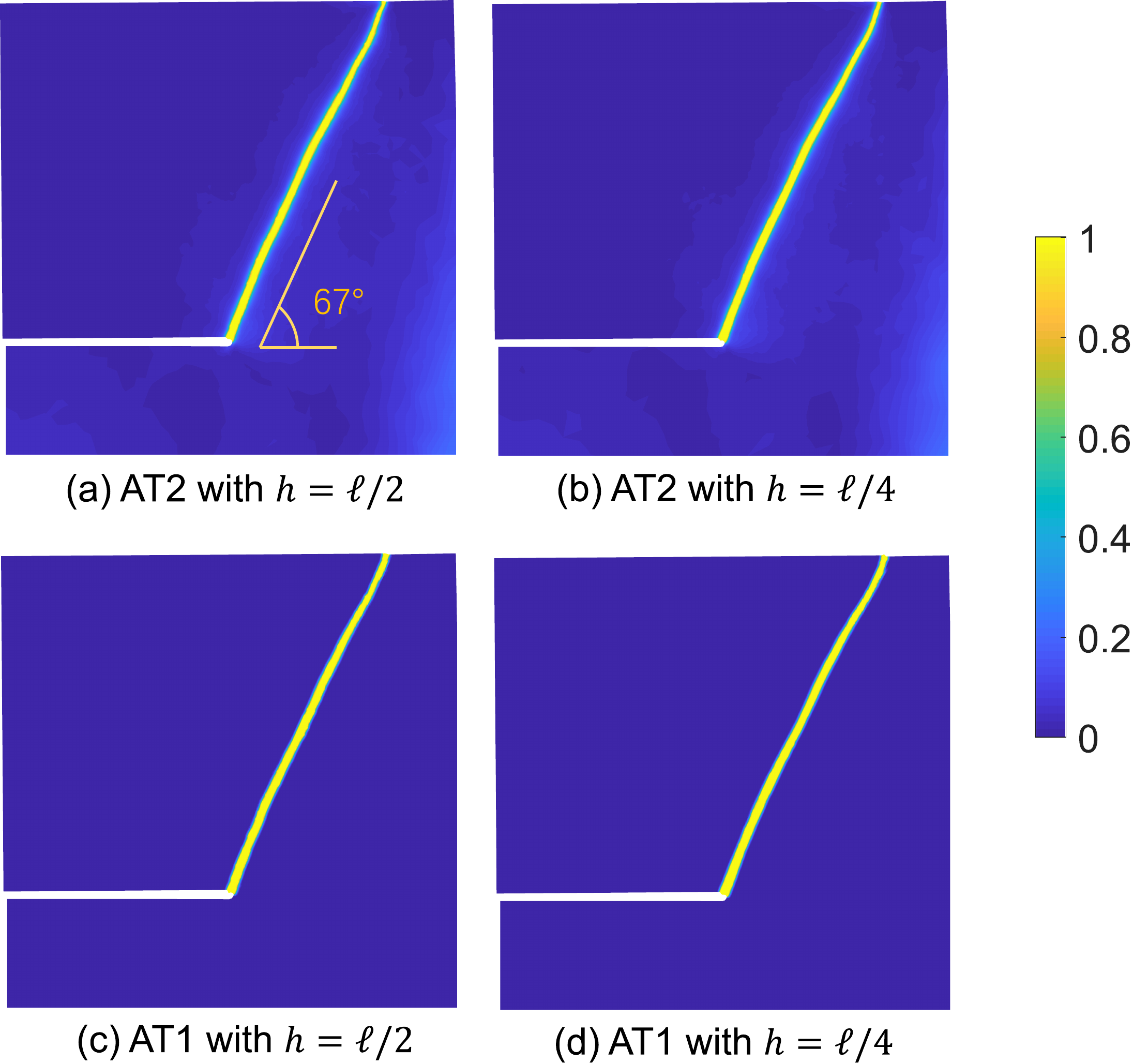}
    \caption{Phase field results for the Kalthoff test of AT1/2 models at time $ t = 87~\mu $s.}
    \label{fig:pf_Kalth}
\end{figure}

Figure \ref{fig:pf_Kalth} shows the final phase field patterns at $t = 87$~$\mu$s for different meshes and models. The crack propagates at 25.5~$\mu$s and with an angle of about $67 ^\circ$ with the horizontal line, which is in good agreement with the experimental results \cite{kalthoff1988failure} and the numerical results using the  phase field method \cite{Chu2017study,Zhou2018phase}.  

\begin{figure}[tbp]
    \centering
    \includegraphics[width=0.95\textwidth]{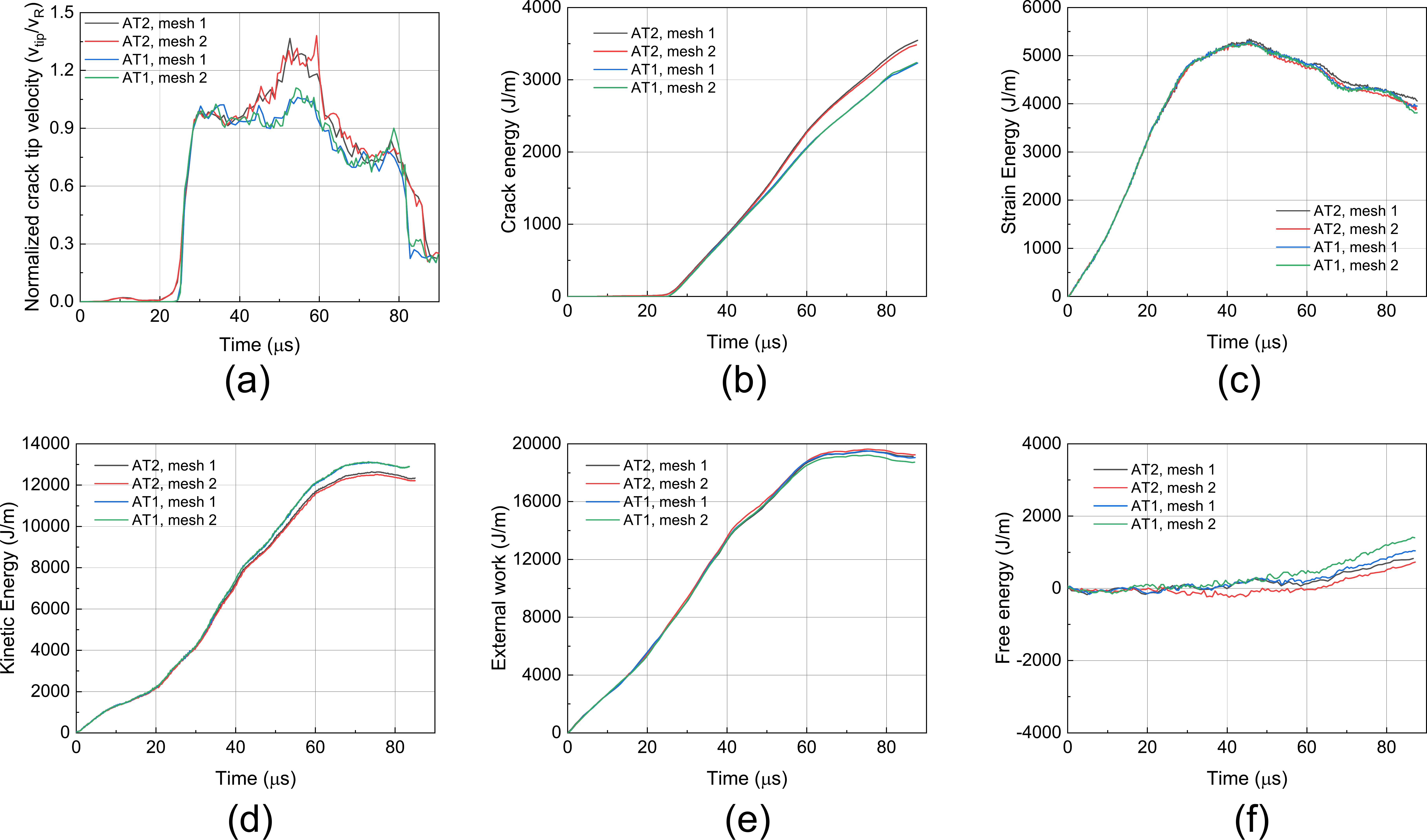}
    \caption{Results of the Kalthoff test. Evolution of (a) normalized crack tip velocity; (b) crack energy; (c) strain energy; (d) kinetic energy; (e) external work; (f) free energy.  The external work (e) is obtained by sampling the power of the reaction force, and then integrating this power with respect to time. The free energy  gradually increases, reaching between 3.74\% and 7.34\% of the external work in the end, which seems to violate the law of conservation of energy. See Figure \ref{fig:gcL} for a comparison.}
    \label{fig:KalthCE}
\end{figure}

Figure \ref{fig:KalthCE}(a) presents the normalized crack tip velocity of the Kalthoff test, the velocity is almost two times the result reported by Liu et al.~\cite{LIU2016Abaqus} The differences in crack tip velocity may be caused by the different post-processing methods, where Ref.~\citenum{LIU2016Abaqus} employed an alternative method that is different from ours by Eq.~\eqref{eq:speed}, to be discussed later. Figure \ref{fig:KalthCE}(b) shows the crack energy calculated by \eqref{eq:Gamma_d}, which agrees well with the numerical results  in Ref.~\citenum{GEELEN2019680}. In addition, the crack energy of the AT2 model is a little higher than the AT1 model for both meshes. Figure \ref{fig:KalthCE}(c) and (d) show the evolution of the strain energy and kinetic energy, respectively, and the results are consistent with numerical results reported by Zhang et al.~\cite{ZHANG2021Quasi} Figure \ref{fig:KalthCE}(e) shows the external work, to the best of our knowledge, there is no relevant report on external work of the Kalthoff test by the phase field method at present. However, our result is in good agreement with the result using the cohesive zone model by Park et al.~\cite{Park2012Adaptive} Figure \ref{fig:KalthCE}(f) shows the evolution of the free energy. As we can see, the free energy gradually increases, reaching between 715.34 J and 1403.51 J at the end of the simulation, which seems to violate the law of conservation of energy. This phenomenon appears to be an open question. 

\begin{figure}[htbp]
    \centering
    \includegraphics[width=0.6\textwidth]{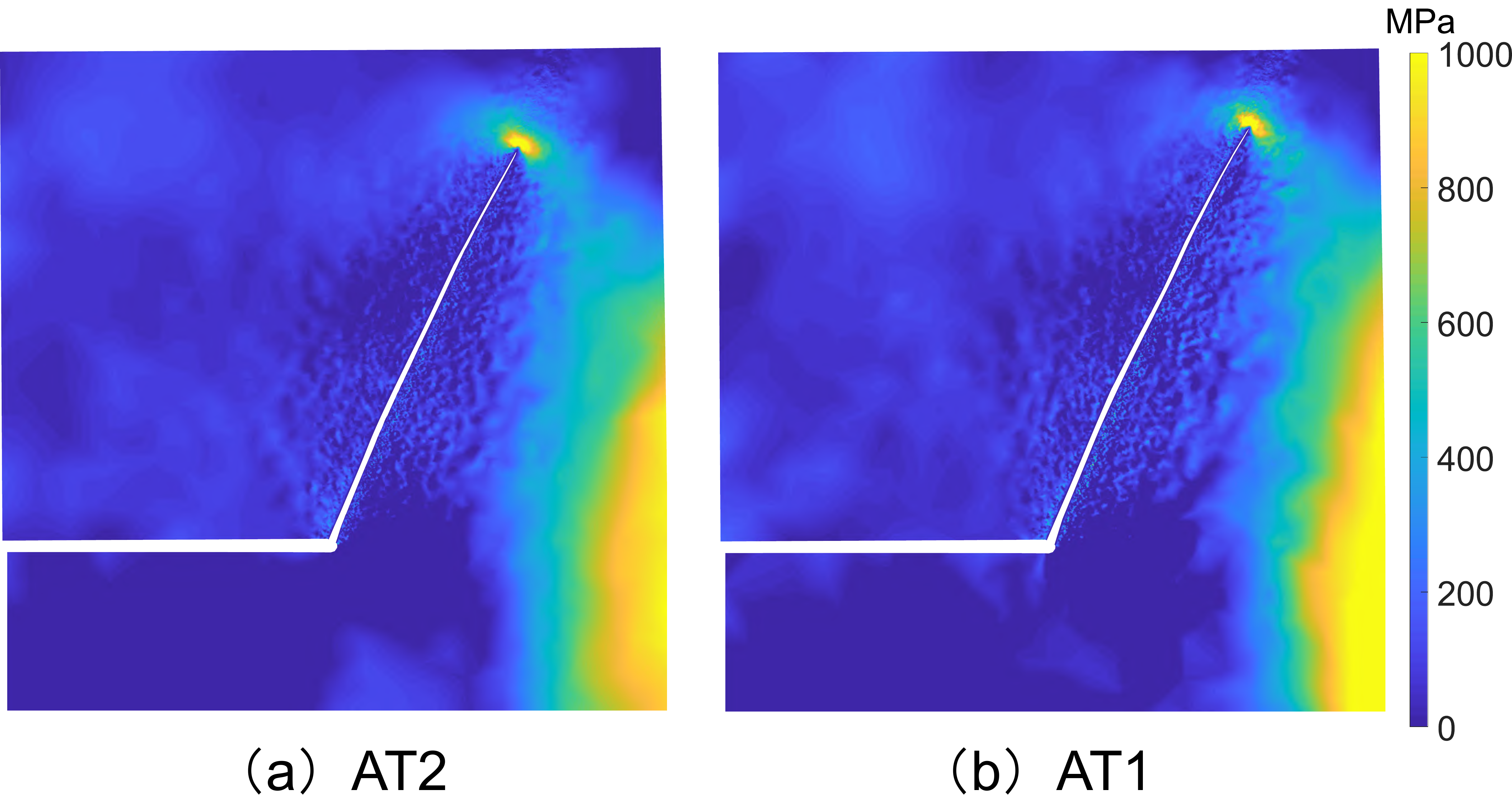}
    \caption{Maximum principal stress of Kalthoff test for Mesh 2 at $t = 70~\mu $s.}
    \label{fig:MaxStr_kalth}
\end{figure}

Figure \ref{fig:MaxStr_kalth} shows the distribution of the maximum principal stress for Mesh 2. The stress concentration is clearly seen at the crack tip, also the bottom right corner, in both AT2 and AT1 models. The result is in good agreement with those in Liu et al.~\cite{LIU2016Abaqus} 

\begin{figure}[htbp]
    \centering
    \includegraphics[width=0.95\textwidth]{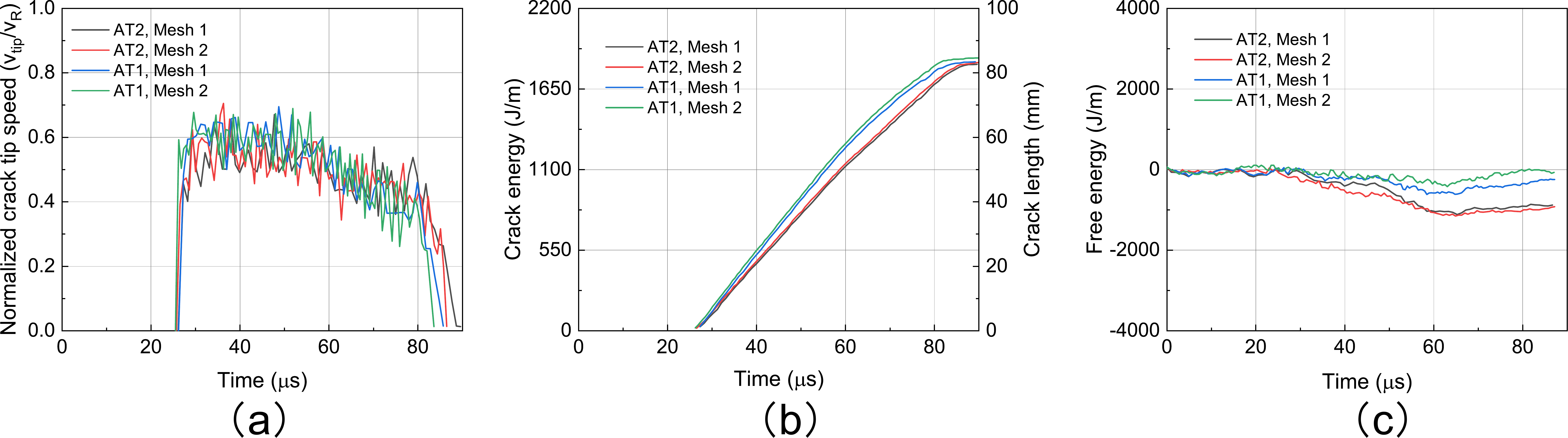}
    \caption{Results of the Kalthoff test. Evolution of (a) normalized crack tip velocity calculated by iso-curve strategy; (b) crack energy calculated by Eq.~\eqref{eq:E_dis}, i.e., $\hat{\Gamma}$ (left vertical axis), and crack length (right vertical axis); (c) free energy $T(\boldsymbol{\dot u}) + V(\boldsymbol{u}, d) + \hat{\Gamma}(d)$, which satisfies the conservation of energy. }
    \label{fig:gcL}
\end{figure}

\paragraph{Alternative method to calculate the crack tip velocity}
As an attempt to reconcile the discrepancy, we employ the iso-curve strategy to calculate the crack tip velocity, as also done by Liu et al.~\cite{LIU2016Abaqus} In this approach, the position of the crack tip is determined by the iso-curve  with phase field $d=0.9$. Therefore, the crack tip velocity is recalculated by $v_n = \|\boldsymbol{x}_{n} - \boldsymbol{x}_{n-1}\| / (t_{n} - t_{n-1})$, where $\boldsymbol{x}_n$ represents the location of the crack tip at $n$th sampling time  $t_n$, and  the result is shown in Figure  \ref{fig:gcL}(a). As can be seen, the crack accelerates to near 0.6$v_R$ and then remains with this velocity during the propagation until it reaches the top boundary,  which agrees well with the result reported  in Ref.~\citenum{LIU2016Abaqus}. 

With this iso-curve scheme, the four cases of Kalthoff test show a similar final crack length of approximately $l_{crack} = 83$ mm, see Figure \ref{fig:gcL}(b) with the right vertical axis. Correspondingly, the crack energy can be computed as
\begin{equation}\label{eq:E_dis}
    \hat{\Gamma} = g_c  l_{crack},
\end{equation}
with the value of 1836.79 J for a sharp crack, see Figure  \ref{fig:gcL}(b) with the left vertical axis. A significant difference is that $\hat{\Gamma}$ is much smaller than $\Gamma $, and the ratios of $\Gamma/\hat{\Gamma}$ are 1.9 and 1.75 for the AT2 and AT1 models, respectively. 

In addition, we recalculate the free energy by using  \eqref{eq:E_dis} instead of \eqref{eq:Gamma_d}, i.e., $T(\boldsymbol{\dot u}) + V(\boldsymbol{u}, d) + \hat{\Gamma}$, and the result is shown in  Figure  \ref{fig:gcL}(c). As seen, with $\hat{\Gamma} $, the results are energetically stable and satisfy the conservation of energy. 

\paragraph{Discussions}
In the Kalthoff test, the crack energy calculated by \eqref{eq:Gamma_d} is higher than that by \eqref{eq:E_dis}. This phenomenon is not unique to this work but also reported in Refs.~\citenum{BORDEN201277,LIU2016Abaqus,Zhou2018phase,GEELEN2019680,ZHANG2021Quasi,Tangella2022Hybrid,Reddy2021Modeling}, in which the ratio of $\Gamma/\hat{\Gamma}$ is between 1.90 and 2.45, equal to or even higher than our value. Meanwhile, in Ref.~\citenum{Vinh2018modeling}, this ratio is 1.37. In addition, this phenomenon was also reported in other dynamic phase field fracture by Ziaei-Rad and Shen \cite{ZIAEIRAD2016Massive}, where the ratio is approximately 2. 

Although the main reason why $\Gamma $ is higher than $\hat{\Gamma}$ need to be further investigated, we suggest that the way of enforcing irreversibility constraint is not an ideal candidate. Borden et al.~\cite{BORDEN201277} suggested that the strain-history field (alternative way to enforce the irresversibility) could play an important role,  but the ratio of $\Gamma/\hat{\Gamma}$ is still 1.4 despite the strain-history field being removed which allows the crack to heal. In addition, Geelen et al.~\cite{GEELEN2019680} employed the augmented Lagrangian method to enforce the irresversibility and the resulting ratio is 2. 

Moreover, Li et al.~\cite{Li2016Gradient} stated that the numerical phase field of the Kalthoff test is wider than the analytical one, and the wider damage profile will lead to an amplified effective fracture toughness, which had also been reported by Bourdin et al.~\cite{Bourdin2008variationl} Furthermore, Bleyer et al.~\cite{Bleyer2017Dynamic} suggested that the mesh size has an influence on the result of both quasi-static and dynamic fracture and that will further lead to an overestimated crack energy (see Eqs.~(16) and (17) in Ref.~\citenum{Bleyer2017Dynamic}  for more details). 


This issue appears to be an open question for the Kalthoff test.


\section{Conclusions}\label{sec5}
In this paper, we have proposed an asynchronous variational formulation for the phase field approach to dynamic brittle fracture. The formulation is derived from Hamilton's principle of stationary action and to a great extent, retains the advantages of variational integrators. A major characteristic of the formulation is that it allows elements to have independent time steps. The result indicates that the  formulation is able to simulate the dynamic fracture propagation and branching successfully.  As a result of the variational structure, the formulation performs a remarkable energy behavior during the simulation. Compared to synchronous methods, the presented formulation improves the computational efficiency for problems involving a high contrast in element sizes or material properties.

Another characteristic is that the phase field irresversibility condition is enforced by the reduced-space active set method at the level of element patches. As a result, AT2 and AT1 variants of the phase field approach may be implemented with similar costs. 
The present study shows that these two variants lead to similar results at roughly the same computational cost.



\section*{Acknowledgments}
We acknowledge the financial support by the National Natural Science Foundation of China, Grant No.~11972227, and by the Natural Science Foundation of Shanghai, Grant No.~19ZR1424200.

\appendix

\section{Phase field result without using element patches} \label{app1}

\begin{figure}[htbp]
    \centering
    \includegraphics[width=0.9\textwidth]{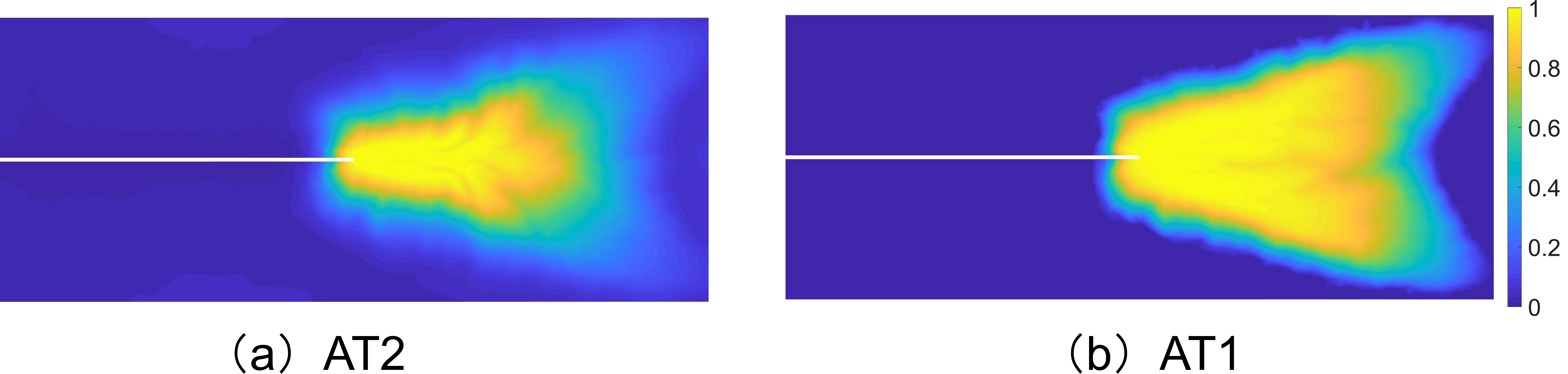}
    \caption{Phase field result with the formulation without using patches.}
    \label{fig:app_unpatch}
\end{figure}

Figure \ref{fig:app_unpatch} shows the phase field result obtained with \eqref{Eq:update_d}, i.e., without patches, the boundary conditions and material properties are the same as those of Section \ref{subsec_branch}. As seen, the crack patterns are too diffused.

\bibliography{wileyNJD-AMA}%

\begin{thebibliography}{10}
\providecommand \doibase [0]{http://dx.doi.org/}%

\bibitem{freund1998dynamic}
Freund LB. {\it Dynamic Fracture Mechanics}.
\newblock Cambridge University Press .
\newblock 1998.

\bibitem{ravi2004dynamic}
Ravi-Chandar K. {\it Dynamic Fracture}.
\newblock Elsevier .
\newblock 2004.

\bibitem{fineberg2015recent}
Fineberg J, Bouchbinder E. Recent developments in dynamic fracture: Some
  perspectives. {\it International Journal of Fracture} 2015\string;
  196(1-2)\string: 33--57.
\newblock \href {\doibase https://doi.org/10.1007/s10704-015-0038-x} {doi:
  https://doi.org/10.1007/s10704-015-0038-x}

\bibitem{SUN2021Astate}
Sun Y, Edwards MG, Chen B, Li C. A state-of-the-art review of crack branching.
  {\it Engineering Fracture Mechanics} 2021\string; 257\string: 108036.
\newblock \href {\doibase https://doi.org/10.1016/j.engfracmech.2021.108036}
  {doi: https://doi.org/10.1016/j.engfracmech.2021.108036}

\bibitem{rabczuk2013computational}
Rabczuk T. Computational methods for fracture in brittle and quasi-brittle
  solids: State-of-the-art review and future perspectives. {\it International
  Scholarly Research Notices} 2013\string; 2013\string: 1-38.
\newblock \href {\doibase https://doi.org/10.1155/2013/849231} {doi:
  https://doi.org/10.1155/2013/849231}

\bibitem{reth2005An}
Réthoré J, Gravouil A, Combescure A. An energy-conserving scheme for dynamic
  crack growth using the eXtended finite element method. {\it International
  Journal for Numerical Methods in Engineering} 2005\string; 63(5)\string:
  631-659.
\newblock \href {\doibase https://doi.org/10.1002/nme.1283} {doi:
  https://doi.org/10.1002/nme.1283}

\bibitem{Nguyen2014Discontinuous}
Nguyen VP. Discontinuous Galerkin/extrinsic cohesive zone modeling:
  Implementation caveats and applications in computational fracture mechanics.
  {\it Engineering Fracture Mechanics} 2014\string; 128\string: 37-68.
\newblock \href {\doibase https://doi.org/10.1016/j.engfracmech.2014.07.003}
  {doi: https://doi.org/10.1016/j.engfracmech.2014.07.003}

\bibitem{Liu2014ARegular}
Liu Y, Filonova V, Hu N, et al. A regularized phenomenological multiscale
  damage model. {\it International Journal for Numerical Methods in
  Engineering} 2014\string; 99(12)\string: 867-887.
\newblock \href {\doibase https://doi.org/10.1002/nme.4705} {doi:
  https://doi.org/10.1002/nme.4705}

\bibitem{ZHANG2019cracking}
Zhang Y, Zhuang X. Cracking elements method for dynamic brittle fracture. {\it
  Theoretical and Applied Fracture Mechanics} 2019\string; 102\string: 1-9.
\newblock \href {\doibase https://doi.org/10.1016/j.tafmec.2018.09.015} {doi:
  https://doi.org/10.1016/j.tafmec.2018.09.015}

\bibitem{song2006Amethod}
Song JH, Areias PMA, Belytschko T. A method for dynamic crack and shear band
  propagation with phantom nodes. {\it International Journal for Numerical
  Methods in Engineering} 2006\string; 67(6)\string: 868-893.
\newblock \href {\doibase https://doi.org/10.1002/nme.1652} {doi:
  https://doi.org/10.1002/nme.1652}

\bibitem{Li2016Gradient}
Li T, Marigo JJ, Guilbaud D, Potapov S. Gradient damage modeling of brittle
  fracture in an explicit dynamics context. {\it International Journal for
  Numerical Methods in Engineering} 2016\string; 108(11)\string: 1381-1405.
\newblock \href {\doibase https://doi.org/10.1002/nme.5262} {doi:
  https://doi.org/10.1002/nme.5262}

\bibitem{Moreau2015Explicit}
Moreau K, Moës N, Picart D, Stainier L. Explicit dynamics with a non-local
  damage model using the thick level set approach. {\it International Journal
  for Numerical Methods in Engineering} 2015\string; 102(3-4)\string: 808-838.
\newblock \href {\doibase https://doi.org/10.1002/nme.4824} {doi:
  https://doi.org/10.1002/nme.4824}

\bibitem{BOURDIN2000797}
Bourdin B, Francfort G, Marigo JJ. Numerical experiments in revisited brittle
  fracture. {\it Journal of the Mechanics and Physics of Solids} 2000\string;
  48(4)\string: 797-826.
\newblock \href {\doibase https://doi.org/10.1016/S0022-5096(99)00028-9} {doi:
  https://doi.org/10.1016/S0022-5096(99)00028-9}

\bibitem{FRANCFORT1998revisit}
Francfort G, Marigo JJ. Revisiting brittle fracture as an energy minimization
  problem. {\it Journal of the Mechanics and Physics of Solids} 1998\string;
  46(8)\string: 1319-1342.
\newblock \href {\doibase https://doi.org/10.1016/S0022-5096(98)00034-9} {doi:
  https://doi.org/10.1016/S0022-5096(98)00034-9}

\bibitem{AMIRI2014thinshell}
Amiri F, Millán D, Shen Y, Rabczuk T, Arroyo M. Phase-field modeling of
  fracture in linear thin shells. {\it Theoretical and Applied Fracture
  Mechanics} 2014\string; 69\string: 102-109.
\newblock \href {\doibase https://doi.org/10.1016/j.tafmec.2013.12.002} {doi:
  https://doi.org/10.1016/j.tafmec.2013.12.002}

\bibitem{LAI2020phase}
Lai W, Gao J, Li Y, Arroyo M, Shen Y. Phase field modeling of brittle fracture
  in an Euler–Bernoulli beam accounting for transverse part-through cracks.
  {\it Computer Methods in Applied Mechanics and Engineering} 2020\string;
  361\string: 112787.
\newblock \href {\doibase https://doi.org/10.1016/j.cma.2019.112787} {doi:
  https://doi.org/10.1016/j.cma.2019.112787}

\bibitem{MOLLAALI2019Numerical}
Mollaali M, Ziaei-Rad V, Shen Y. Numerical modeling of CO2 fracturing by the
  phase field approach. {\it Journal of Natural Gas Science and Engineering}
  2019\string; 70\string: 102905.
\newblock \href {\doibase https://doi.org/10.1016/j.jngse.2019.102905} {doi:
  https://doi.org/10.1016/j.jngse.2019.102905}

\bibitem{shen2018implementation}
Shen Y, Mollaali M, Li Y, Ma W, Jiang J. Implementation details for the phase
  field approaches to fracture. {\it Journal of Shanghai Jiaotong University
  (Science)} 2018\string; 23(1)\string: 166--174.
\newblock \href {\doibase https://doi.org/10.1007/s12204-018-1922-0} {doi:
  https://doi.org/10.1007/s12204-018-1922-0}

\bibitem{BORDEN201277}
Borden MJ, Verhoosel CV, Scott MA, Hughes TJ, Landis CM. A phase-field
  description of dynamic brittle fracture. {\it Computer Methods in Applied
  Mechanics and Engineering} 2012\string; 217-220\string: 77-95.
\newblock \href {\doibase https://doi.org/10.1016/j.cma.2012.01.008} {doi:
  https://doi.org/10.1016/j.cma.2012.01.008}

\bibitem{Vinh2018modeling}
Nguyen VP, Wu JY. Modeling dynamic fracture of solids with a phase-field
  regularized cohesive zone model. {\it Computer Methods in Applied Mechanics
  and Engineering} 2018\string; 340\string: 1000-1022.
\newblock \href {\doibase https://doi.org/10.1016/j.cma.2018.06.015} {doi:
  https://doi.org/10.1016/j.cma.2018.06.015}

\bibitem{HAO2022phasefield}
Hao S, Shen Y, Cheng JB. Phase field formulation for the fracture of a metal
  under impact with a fluid formulation. {\it Engineering Fracture Mechanics}
  2022\string; 261\string: 108142.
\newblock \href {\doibase https://doi.org/10.1016/j.engfracmech.2021.108142}
  {doi: https://doi.org/10.1016/j.engfracmech.2021.108142}

\bibitem{HAO2022aphase}
Hao S, Chen Y, Cheng JB, Shen Y. A phase field model for high-speed impact
  based on the updated Lagrangian formulation. {\it Finite Elements in Analysis
  and Design} 2022\string; 199\string: 103652.
\newblock \href {\doibase https://doi.org/10.1016/j.finel.2021.103652} {doi:
  https://doi.org/10.1016/j.finel.2021.103652}

\bibitem{Tian2019hybird}
Tian F, Tang X, Xu T, Yang J, Li L. A hybrid adaptive finite element
  phase-field method for quasi-static and dynamic brittle fracture. {\it
  International Journal for Numerical Methods in Engineering} 2019\string;
  120(9)\string: 1108-1125.
\newblock \href {\doibase https://doi.org/10.1002/nme.6172} {doi:
  https://doi.org/10.1002/nme.6172}

\bibitem{ZIAEIRAD2016Massive}
Ziaei-Rad V, Shen Y. Massive parallelization of the phase field formulation for
  crack propagation with time adaptivity. {\it Computer Methods in Applied
  Mechanics and Engineering} 2016\string; 312\string: 224-253.
\newblock \href {\doibase https://doi.org/10.1016/j.cma.2016.04.013} {doi:
  https://doi.org/10.1016/j.cma.2016.04.013}

\bibitem{Li2019variational}
Li Y, Lai W, Shen Y. Variational h-adaption method for the phase field approach
  to fracture. {\it International Journal of Fracture} 2019\string; 217\string:
  83-103.
\newblock \href {\doibase https://doi.org/10.1007/s10704-019-00372-y} {doi:
  https://doi.org/10.1007/s10704-019-00372-y}

\bibitem{Thomas2021Lscheme}
Engwer C, Pop IS, Wick T. Dynamic and weighted stabilizations of the L-scheme
  applied to a phase-field model for fracture propagation. In: Springer
  International Publishing; 2021; Cham\string: 1177--1184

\bibitem{lew2003asynchronous}
Lew A, Marsden JE, Ortiz M, West M. Asynchronous variational integrators. {\it
  Archive for Rational Mechanics and Analysis} 2003\string; 167(2)\string:
  85--146.
\newblock \href {\doibase https://doi.org/10.1007/s00205-002-0212-y} {doi:
  https://doi.org/10.1007/s00205-002-0212-y}

\bibitem{lew2004variational}
Lew A, Marsden J, Ortiz M, West M. Variational time integrators. {\it
  International Journal for Numerical Methods in Engineering} 2004\string;
  60(1)\string: 153--212.
\newblock \href {\doibase https://doi.org/10.1002/nme.958} {doi:
  https://doi.org/10.1002/nme.958}

\bibitem{FONG2008Stability}
Fong W, Darve E, Lew A. Stability of asynchronous variational integrators. {\it
  Journal of Computational Physics} 2008\string; 227(18)\string: 8367-8394.
\newblock \href {\doibase https://doi.org/10.1016/j.jcp.2008.05.017} {doi:
  https://doi.org/10.1016/j.jcp.2008.05.017}

\bibitem{Matteo2008Convergence}
Focardi M, Mariano PM. Convergence of asynchronous variational integrators in
  linear elastodynamics. {\it International Journal for Numerical Methods in
  Engineering} 2008\string; 75(7)\string: 755-769.
\newblock \href {\doibase https://doi.org/10.1002/nme.2271} {doi:
  https://doi.org/10.1002/nme.2271}

\bibitem{Ryckman2012AVIcontact}
Ryckman RA, Lew AJ. An explicit asynchronous contact algorithm for elastic
  body-rigid wall interaction. {\it International Journal for Numerical Methods
  in Engineering} 2012\string; 89(7)\string: 869-896.
\newblock \href {\doibase https://doi.org/10.1002/nme.3266} {doi:
  https://doi.org/10.1002/nme.3266}

\bibitem{Liu2020Extended}
Liu P, Yang JZ, Yuan C. Extended synchronous variational integrators for wave
  propagations on non-uniform meshes. {\it Communications in Computational
  Physics} 2020\string; 28(2)\string: 691--722.
\newblock \href {\doibase https://doi.org/10.4208/cicp.OA-2019-0167} {doi:
  https://doi.org/10.4208/cicp.OA-2019-0167}

\bibitem{Thomas2008asyncloth}
Thomaszewski B, Pabst S, Stra{\ss}er W. Asynchronous Cloth Simulation. In: The
  2008 International Conference on Computer Graphics and Virtual Reality; 2008.

\bibitem{Pham2011Gradient}
Pham K, Amor H, Marigo JJ, Maurini C. Gradient damage models and their use to
  approximate brittle fracture. {\it International Journal of Damage Mechanics}
  2011\string; 20(4)\string: 618-652.
\newblock \href {\doibase 10.1177/1056789510386852} {doi:
  10.1177/1056789510386852}

\bibitem{REN2019explict}
Ren H, Zhuang X, Anitescu C, Rabczuk T. An explicit phase field method for
  brittle dynamic fracture. {\it Computers \& Structures} 2019\string;
  217\string: 45-56.
\newblock \href {\doibase https://doi.org/10.1016/j.compstruc.2019.03.005}
  {doi: https://doi.org/10.1016/j.compstruc.2019.03.005}

\bibitem{SUH2021asynchronous}
Suh HS, Sun W. Asynchronous phase field fracture model for porous media with
  thermally non-equilibrated constituents. {\it Computer Methods in Applied
  Mechanics and Engineering} 2021\string; 387\string: 114182.
\newblock \href {\doibase https://doi.org/10.1016/j.cma.2021.114182} {doi:
  https://doi.org/10.1016/j.cma.2021.114182}

\bibitem{MIEHE20102765}
Miehe C, Hofacker M, Welschinger F. A phase field model for rate-independent
  crack propagation: Robust algorithmic implementation based on operator
  splits. {\it Computer Methods in Applied Mechanics and Engineering}
  2010\string; 199(45)\string: 2765-2778.
\newblock \href {\doibase 10.1016/j.cma.2010.04.011} {doi:
  10.1016/j.cma.2010.04.011}

\bibitem{AMOR20091209}
Amor H, Marigo JJ, Maurini C. Regularized formulation of the variational
  brittle fracture with unilateral contact: Numerical experiments. {\it Journal
  of the Mechanics and Physics of Solids} 2009\string; 57(8)\string: 1209-1229.
\newblock \href {\doibase https://doi.org/10.1016/j.jmps.2009.04.011} {doi:
  https://doi.org/10.1016/j.jmps.2009.04.011}

\bibitem{LIU2021107358}
Liu Y, Cheng C, Ziaei-Rad V, Shen Y. A micromechanics-informed phase field
  model for brittle fracture accounting for unilateral constraint. {\it
  Engineering Fracture Mechanics} 2021\string; 241\string: 107358.
\newblock \href {\doibase https://doi.org/10.1016/j.engfracmech.2020.107358}
  {doi: https://doi.org/10.1016/j.engfracmech.2020.107358}

\bibitem{WU2020112629}
Wu JY, Nguyen VP, Zhou H, Huang Y. A variationally consistent phase-field
  anisotropic damage model for fracture. {\it Computer Methods in Applied
  Mechanics and Engineering} 2020\string; 358\string: 112629.
\newblock \href {\doibase https://doi.org/10.1016/j.cma.2019.112629} {doi:
  https://doi.org/10.1016/j.cma.2019.112629}

\bibitem{marsden_2001}
Marsden JE, West M. Discrete mechanics and variational integrators. {\it Acta
  Numerica} 2001\string; 10(10)\string: 357-514.
\newblock \href {\doibase 10.1017/S096249290100006X} {doi:
  10.1017/S096249290100006X}

\bibitem{Gerasi2019Penalization}
Gerasimov T, {De Lorenzis} L. On penalization in variational phase-field models
  of brittle fracture. {\it Computer Methods in Applied Mechanics and
  Engineering} 2019\string; 354\string: 990-1026.
\newblock \href {\doibase https://doi.org/10.1016/j.cma.2019.05.038} {doi:
  https://doi.org/10.1016/j.cma.2019.05.038}

\bibitem{GEELEN2019680}
Geelen RJ, Liu Y, Hu T, Tupek MR, Dolbow JE. A phase-field formulation for
  dynamic cohesive fracture. {\it Computer Methods in Applied Mechanics and
  Engineering} 2019\string; 348\string: 680-711.
\newblock \href {\doibase https://doi.org/10.1016/j.cma.2019.01.026} {doi:
  https://doi.org/10.1016/j.cma.2019.01.026}

\bibitem{WU2018Robust}
Wu JY. Robust numerical implementation of non-standard phase-field damage
  models for failure in solids. {\it Computer Methods in Applied Mechanics and
  Engineering} 2018\string; 340\string: 767-797.
\newblock \href {\doibase https://doi.org/10.1016/j.cma.2018.06.007} {doi:
  https://doi.org/10.1016/j.cma.2018.06.007}

\bibitem{Farrell2017Linear}
Farrell P, Maurini C. Linear and nonlinear solvers for variational phase-field
  models of brittle fracture. {\it International Journal for Numerical Methods
  in Engineering} 2017\string; 109(5)\string: 648-667.
\newblock \href {\doibase https://doi.org/10.1002/nme.5300} {doi:
  https://doi.org/10.1002/nme.5300}

\bibitem{Knuth1998Art}
Knuth DE. {\it The Art of Computer Programming}.
\newblock Addison Wesley .
\newblock 1998.

\bibitem{Song2008_comparative}
Song JH, Wang H, Belytschko T. A comparative study on finite element methods
  for dynamic fracture. {\it Computational Mechanics} 2008\string; 42\string:
  239-250.
\newblock \href {\doibase 10.1007/s00466-007-0210-x} {doi:
  10.1007/s00466-007-0210-x}

\bibitem{Ramulu1985Mechanics}
Ramulu M, Kobayashi AS. Mechanics of crack curving and branching -- a dynamic
  fracture analysis. {\it International Journal of Fracture} 1985\string;
  27\string: 187-201.
\newblock \href {\doibase https://doi.org/10.1007/BF00017967} {doi:
  https://doi.org/10.1007/BF00017967}

\bibitem{Sharon1996Micro}
Sharon E, Fineberg J. Microbranching instability and the dynamic fracture of
  brittle materials. {\it Physical Review B Condensed Matter} 1996\string;
  54(10)\string: 7128-7139.
\newblock \href {\doibase 10.1103/physrevb.54.7128} {doi:
  10.1103/physrevb.54.7128}

\bibitem{Schluter2014Phase}
Schlüter A, Kuhn C, Müller R. Phase field approximation of dynamic brittle
  fracture. {\it Proceedings in Applied Mathematics and Mechanics} 2014\string;
  14(1)\string: 143-144.
\newblock \href {\doibase https://doi.org/10.1002/pamm.201410059} {doi:
  https://doi.org/10.1002/pamm.201410059}

\bibitem{LIU2016Abaqus}
Liu G, Li Q, Msekh MA, Zuo Z. Abaqus implementation of monolithic and staggered
  schemes for quasi-static and dynamic fracture phase-field model. {\it
  Computational Materials Science} 2016\string; 121\string: 35-47.
\newblock \href {\doibase https://doi.org/10.1016/j.commatsci.2016.04.009}
  {doi: https://doi.org/10.1016/j.commatsci.2016.04.009}

\bibitem{Bobaru2015Why}
Bobaru F, Zhang G. Why do cracks branch? A peridynamic investigation of dynamic
  brittle fracture. {\it International Journal of Fracture} 2015\string;
  196\string: 59-98.
\newblock \href {\doibase 10.1007/s10704-015-0056-8} {doi:
  10.1007/s10704-015-0056-8}

\bibitem{MANDAL2020evaluation}
Mandal TK, Nguyen VP, Wu JY. Evaluation of variational phase-field models for
  dynamic brittle fracture. {\it Engineering Fracture Mechanics} 2020\string;
  235\string: 107169.
\newblock \href {\doibase https://doi.org/10.1016/j.engfracmech.2020.107169}
  {doi: https://doi.org/10.1016/j.engfracmech.2020.107169}

\bibitem{kalthoff1988failure}
Kalthoff JF. Shadow optical analysis of dynamic shear fracture. {\it Optical
  Engineering} 1988\string; 27(10)\string: 835 -- 840.
\newblock \href {\doibase 10.1117/12.7976772} {doi: 10.1117/12.7976772}

\bibitem{kalthoff2000modes}
Kalthoff JF. Modes of dynamic shear failure in solids. {\it International
  Journal of Fracture} 2000\string; 101(1)\string: 1--31.
\newblock \href {\doibase https://doi.org/10.1023/A:1007647800529} {doi:
  https://doi.org/10.1023/A:1007647800529}

\bibitem{WANG2019nonordinary}
Wang H, Xu Y, Huang D. A non-ordinary state-based peridynamic formulation for
  thermo-visco-plastic deformation and impact fracture. {\it International
  Journal of Mechanical Sciences} 2019\string; 159\string: 336-344.
\newblock \href {\doibase https://doi.org/10.1016/j.ijmecsci.2019.06.008} {doi:
  https://doi.org/10.1016/j.ijmecsci.2019.06.008}

\bibitem{Chu2017study}
Chu D, Li X, Liu Z. Study the dynamic crack path in brittle material under
  thermal shock loading by phase field modeling. {\it International Journal of
  Fracture} 2017\string; 208\string: 115-130.
\newblock \href {\doibase https://doi.org/10.1007/s10704-017-0220-4} {doi:
  https://doi.org/10.1007/s10704-017-0220-4}

\bibitem{Zhou2018phase}
Zhou S, Rabczuk T, Zhuang X. Phase field modeling of quasi-static and dynamic
  crack propagation: COMSOL implementation and case studies. {\it Advances in
  Engineering Software} 2018\string; 122\string: 31-49.
\newblock \href {\doibase https://doi.org/10.1016/j.advengsoft.2018.03.012}
  {doi: https://doi.org/10.1016/j.advengsoft.2018.03.012}

\bibitem{ZHANG2021Quasi}
Zhang Y, Ren H, Areias P, Zhuang X, Rabczuk T. Quasi-static and dynamic
  fracture modeling by the nonlocal operator method. {\it Engineering Analysis
  with Boundary Elements} 2021\string; 133\string: 120-137.
\newblock \href {\doibase https://doi.org/10.1016/j.enganabound.2021.08.020}
  {doi: https://doi.org/10.1016/j.enganabound.2021.08.020}

\bibitem{Park2012Adaptive}
Park K, Paulino GH, Celes W, Espinha R. Adaptive mesh refinement and coarsening
  for cohesive zone modeling of dynamic fracture. {\it International Journal
  for Numerical Methods in Engineering} 2012\string; 92(1)\string: 1-35.
\newblock \href {\doibase https://doi.org/10.1002/nme.3163} {doi:
  https://doi.org/10.1002/nme.3163}

\bibitem{Tangella2022Hybrid}
Tangella RG, Kumbhar P, Annabattula RK. Hybrid Phase Field Modelling of Dynamic
  Brittle Fracture and Implementation in FEniCS. In:  Krishnapillai S, R. V, Ha
  SK. \kern-2pt, eds. {\it Composite Materials for Extreme Loading}Springer
  Singapore; 2022; Singapore\string: 15--24

\bibitem{Reddy2021Modeling}
Reddy SSK, Amirtham R, Reddy JN. Modeling fracture in brittle materials with
  inertia effects using the phase field method. {\it Mechanics of Advanced
  Materials and Structures} 2021\string; 0(0)\string: 1-16.
\newblock \href {\doibase 10.1080/15376494.2021.2010289} {doi:
  10.1080/15376494.2021.2010289}

\bibitem{Bourdin2008variationl}
Bourdin B, Francfort GA, Marigo JJ. The variational approach to fracture. {\it
  Journal of Elasticity} 2008\string; 91\string: 5-148.
\newblock \href {\doibase https://doi.org/10.1007/s10659-007-9107-3} {doi:
  https://doi.org/10.1007/s10659-007-9107-3}

\bibitem{Bleyer2017Dynamic}
Bleyer J, Roux-Langlois C, Molinari JF. Dynamic crack propagation with a
  variational phase-field model: Limiting speed, crack branching and
  velocity-toughening mechanisms. {\it International Journal of Fracture}
  2017\string; 204\string: 79-100.
\newblock \href {\doibase https://doi.org/10.1007/s10704-016-0163-1} {doi:
  https://doi.org/10.1007/s10704-016-0163-1}

\end{thebibliography}







\end{document}